\documentclass[letterpaper]{article}
\usepackage[T1]{fontenc}
\usepackage[latin1]{inputenc}
\usepackage{amsmath,amsthm,amsfonts,amssymb,euscript,stmaryrd,graphicx,appendix}
\usepackage{algorithm,algorithmic}
\usepackage{color}
\usepackage{listings}
\usepackage{pstricks,psfrag}
\usepackage{amsthm}
\usepackage{fullpage}
\usepackage{yhmath} 
\usepackage{enumitem}

\definecolor{colKeys}{rgb}{0,0,1}
\definecolor{unidentified}{rgb}{0,0,0}
\definecolor{concealments}{rgb}{0,0.5,0}
\definecolor{colString}{rgb}{0.6,0.1,0.1}

\definecolor{Grey1}{rgb}{0.5,0.5,0.5}
\definecolor{Grey2}{rgb}{0.8,0.8,0.8}
\definecolor{Grey3}{rgb}{0.9,0.9,0.9}
\definecolor{Grey4}{rgb}{0.95,0.95,0.95}
\definecolor{ACORRIGER}{rgb}{0.9,0,0}

\def\myqed {{
\parfillskip=0pt 
\widowpenalty=10000 
\displaywidowpenalty=10000 
\finalhyphendemerits=0 
\leavevmode 
\unskip 
\nobreak 
\hfil 
\penalty50 
\hskip.2em 
\null 
\hfill 
$\square$
\par}} 

\newcommand{\footnoteremember}[2]{
 \footnote{#2}
 \newcounter{#1}
 \setcounter{#1}{\value{footnote}}
}

\def\independent{{\perp\!\!\!\!\perp}}
\def\simdist{\stackrel{\mathcal{L}}{\sim}}

\author
{
Sylvain \textsc{Corlay}\footnote{Bloomberg L.P. Quantitative Finance Research, 731 Lexington avenue, New York, NY 10022, USA. E-mail: sylvain.corlay@gmail.com.}
\footnoteremember{myfootnote}{Laboratoire de Probabilités et Modèles Aléatoires, UMR 7599, Université Paris 6, case 188, 4, pl. Jussieu, F-75252 Paris Cedex 5, France. }
}

\title{Partial functional quantization and generalized bridges}
\date{September 17, 2012}
\newtheorem{theo}{Theorem}[section]
\newtheorem{prop}[theo]{Proposition}
\newtheorem{coro}[theo]{Corollary}
\newtheorem*{remark}{Remark}
\newtheorem*{example}{Example}
\newtheorem{lemm}[theo]{Lemma}

\newcommand{\1}{\textbf{1}}
\newcommand{\N}{\mathbb{N}}

\newcommand{\R}{\mathbb{R}}

\newcommand{\E}{\mathbb{E}}

\newcommand{\PP}{\mathbb{P}}

\newcommand{\card}{\operatorname{card}}

\newcommand{\sspan}{\operatorname{span}}
\newcommand{\cov}{\operatorname{cov}}
\newcommand{\var}{\operatorname{Var}}
\newcommand{\Proj}{\operatorname{Proj}}

\def\keywordname{{\bf Keywords:}} 
\newcommand{\keywords}[1]{\par\addvspace\baselineskip\noindent\keywordname\enspace\ignorespaces#1}

\graphicspath{{./}{figures/}}
\begin{document}
\lstset{language=c++}
\maketitle
\begin{abstract}
\par In this article, we develop a new approach to functional quantization, which consists in discretizing only a finite subset of the Karhunen-Loève coordinates of a continuous Gaussian semimartingale $X$. 
\par Using filtration enlargement techniques, we prove that the conditional distribution of $X$ knowing its first Karhunen-Loève coordinates is a Gaussian semimartingale with respect to a bigger filtration. This allows us to define the partial quantization of a solution of a stochastic differential equation with respect to $X$ by simply plugging the partial functional quantization of $X$ in the SDE. 
\par Then we provide an upper bound of the $L^p$-partial quantization error for the solution of SDEs involving the $L^{p+\varepsilon}$-partial quantization error for $X$, for $\varepsilon >0$. The $a.s.$ convergence is also investigated. 
\par Incidentally, we show that the conditional distribution of a Gaussian semimartingale $X$, knowing that it stands in some given Voronoi cell of its functional quantization, is a (non-Gaussian) semimartingale. As a consequence, the functional stratification method developed in \cite{CorlayPagesStratification} amounted, in the case of solutions of SDEs, to using the Euler scheme of these SDEs in each Voronoi cell.
\end{abstract}
\section*{Introduction}\label{sec:introduction_quantization}
\par \noindent Let $(\Omega,\mathcal{A},\PP)$ be a probability space, and $E$ a reflexive separable Banach space. The norm on $E$ is denoted by $|\cdot |$. The quantization of a $E$-valued random variable $X$ consists in its approximation by a random variable $Y$ taking finitely many values. The resulting error of this discretization is measured by the $L^p$ norm of $|X-Y|$. If we settle on a fixed maximum cardinal for $Y(\Omega)$, the minimization of the quantization error amounts to the minimization problem:
\begin{equation}\label{eq:minimization_quantization_partial}
\min \left\{ \big\| |X-Y| \big\|_p, \ Y: \Omega \to E \textnormal{ measurable}, \ \card(Y(\Omega)) \leq N \right\}.
\end{equation}
\par \noindent A solution to (\ref{eq:minimization_quantization_partial}) is an optimal quantizer of $X$. The corresponding quantization error is denoted by $\mathcal{E}_{N,p}(X) := \min \left\{ \big\| |X-Y| \big\|_p, \ Y: \Omega \to E \textnormal{ measurable}, \ \card(Y(\Omega)) \leq N \right\}$. One usually drops the $p$ subscript in the quadratic case ($p = 2$). This problem, initially investigated as a signal discretization method \cite{GershoGrayVectorQuantization}, has then been introduced in numerical probability to devise cubature methods \cite{PagesIntegVectorQuant} or to solve multidimensional stochastic control problems \cite{BallyPagesPrintemsAmerican1}. Since the early $2000$'s, the infinite-dimensional setting has been extensively investigated from both constructive numerical and theoretical viewpoints with a special attention paid to functional quantization, especially in the quadratic case \cite{LuschgyPagesFunctional3} but also in some other Banach spaces \cite{WilbertzPHD}. Stochastic 
processes are viewed as random variables taking values in functional spaces.
\par We now assume that $X$ is a bi-measurable stochastic process on $[0,T]$ verifying $\int_0^T \E\left[|X_t|^2\right]dt < +\infty$, so that this can be viewed as a random variable valued in the separable Hilbert space $L^2([0,T])$. We assume that its covariance function $\Gamma^X$ is continuous. In the seminal article on Gaussian functional quantization \cite{LuschgyPagesFunctional3}, it is shown that in the centered Gaussian case, linear subspaces $U$ of $L^2([0,T])$ spanned by $L^2$-optimal quantizers correspond to principal components of $X$. In other words, they are spanned by the first eigenvectors of the covariance operator of $X$. Thus, the quadratic optimal quantization of $X$ involves its Karhunen-Loève eigensystem $(e_n^X,\lambda_n^X)_{n \geq 1}$. If $Y$ is a quadratic $N$-optimal quantizer of $X$ and $d^X(N)$ is the dimension of the subspace of $L^2([0,T])$ spanned by $Y(\Omega)$, the quadratic quantization error $\mathcal{E}_N^2(X)$ verifies
\begin{equation}\label{eq:distorsion_representation}
\mathcal{E}_N^2(X) = \sum\limits_{j \geq m+1} \lambda_j^X + \mathcal{E}_N^2\left( \bigotimes\limits_{j=1}^m \mathcal{N}\left(0,\lambda_j^X\right)\right) \textnormal{ for } m \geq d^X(N).
\end{equation}
\begin{equation}\label{eq:distorsion_representation_2}
\mathcal{E}_N^2(X) < \sum\limits_{j \geq m+1} \lambda_j^X + \mathcal{E}_N^2\left( \bigotimes\limits_{j=1}^m \mathcal{N}\left(0,\lambda_j^X\right)\right) \textnormal{ for } 1 \leq m < d^X(N).
\end{equation}
\par
\begin{figure}[!ht]
	\begin{center}
	\psfrag{-2}{$-2$}
	\psfrag{-1.5}{$-1.5$}
	\psfrag{-1}{$-1$}
	\psfrag{-0.5}{$-0.5$}
	\psfrag{0}{$0$}
	\psfrag{0.5}{$0.5$}
	\psfrag{1}{$1$}
	\psfrag{1.5}{$1.5$}
	\psfrag{2}{$2$}
	\psfrag{0.2}{$0.2$}
	\psfrag{0.4}{$0.4$}
	\psfrag{0.6}{$0.6$}
	\psfrag{0.8}{$0.8$}
	\includegraphics[height=6cm]{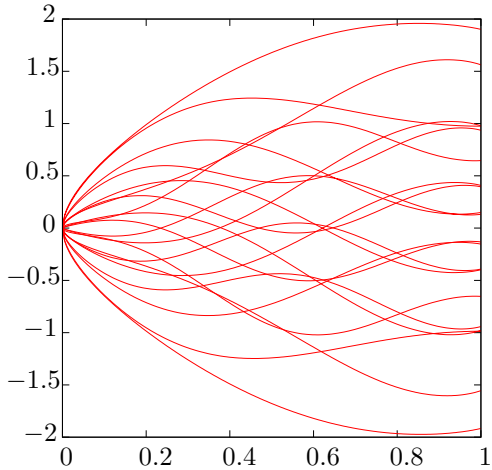}
	\caption{Quadratic $N$-optimal quantizer of the fractional Brownian motion on $[0,1]$ with Hurst parameter $H = 0.25$ and $N=20$. The quantization dimension is $3$.}
	\label{fig:fractional_brownian_motion_optimal_quantization_partial}
	\end{center}
\end{figure}

 \noindent To perform optimal quantization, the decomposition is first truncated at a fixed order $m$ and then the $\R^m$-valued Gaussian vector, constituted of the $m$ first coordinates of the process on its Karhunen-Loève decomposition, is quantized. To reach optimality, we have to determine the optimal rank of truncation $d^X(N)$ (the quantization dimension) and the optimal $d^X(N)$-dimensional quantizer corresponding to the first coordinates $\bigotimes\limits_{j=1}^{d^X(N)} \mathcal{N}\left(0,\lambda_j^X\right)$. A sharply optimized database of quantizers of univariate and multivariate Gaussian distributions is available on the web site {\verb www.quantize.maths-fi.com } \cite{WebSiteGaussian} for download. Usual examples of such processes are the standard Brownian motion on $[0,T]$, the Brownian bridge on $[0,T]$, Ornstein-Uhlenbeck processes and the fractional Brownian motion. In Figure \ref{fig:fractional_brownian_motion_optimal_quantization_partial}, we display the quadratic optimal $N$-
quantizer of the fractional Brownian motion on $[0,1]$ with Hurst exponent $H = 0.25$ and $N = 20$. 
\par From a constructive viewpoint, the numerical computation of the optimal quantization or the optimal product quantization requires a numerical evaluation of the Karhunen-Loève eigenfunctions and eigenvalues, at least the very first terms. (As seen in \cite{LuschgyPagesFunctional3,LuschgyPagesFunctional2,LuschgyPagesWilberzFunctionalSchemes}, under rather general conditions on its eigenvalues, the quantization dimension of a Gaussian process increases asymptotically as the logarithm of the size of the quantizer. Hence it is most likely that it is small. For instance, the quantization dimension of Brownian motion with $N=10000$ is $9$.) The Karhunen-Loève decompositions of several usual Gaussian processes have a closed-form expression. This is the case for standard Brownian motion, Brownian bridge and Ornstein-Uhlenbeck processes. The case of Ornstein-Uhlenbeck processes is derived in \cite{CorlayPagesStratification}, in the general setting of an arbitrary initial variance $\sigma_0$. Another example of 
explicit Karhunen-Loève expansion is derived in \cite{DeheuvelsMartynov} by Deheuvels and Martynov. 
\par In the general case, no closed-form expression is available for the Karhunen-Loève expansion. For example, the K-L expansion of the fractional Brownian motion is not known. Yet, one can use numerical schemes to solve the correspnding eigenvalue problem. In \cite{CorlayNystrom}, the so-called ``Nyström method'' is used to compute the first terms of the K-L decomposition of the fractional Brownian motion and to perform its optimal quantization. 
\vspace{1mm}
\par In this article, we propose a new functional quantization scheme for a bi-measurable Gaussian process $X$, which consists in discretizing a finite subset of its Karhunen-Loève coordinates, instead of performing a full quantization. This \emph{partial functional quantization} approach is motivated by two observations. The first one is that the conditional distribution of $X$ knowing that it falls into a given $L^2$ Voronoi cell of its optimal quantizer is the crux of the recently developed functional stratification scheme \cite{CorlayPagesStratification}. It comes to conditioning the process with respect to its first Karhunen-Loève coordinates. This work provides a better justification of the functional stratification scheme of \cite{CorlayPagesStratification}. The second observation is that one of the main purposes of the (full) functional quantization of $X$ is to perform a quantization of the solution of a SDE with respect to $X$, when a stochastic integration with respect to $X$ can be defined (see \cite{PagesPrintemsFunctional4, LuschgyPagesFunctional5, PagesSellamiSDE}). As (full) functional quantizers of $X$ will typically have bounded variations, one needs to add a correction term to the SDE. Eventually, this comes to plug the functional quantizer of $X$ in the SDE written in the Stratonovich sense. In contrast, the partial quantization of $X$ can be directly plugged into the SDE written in the Itô sense. We provide $a.s.$ and $L^p$ convergence results for this method.
\vspace{1mm}
\par The paper is organized as follows: Section \ref{seq:quantization_based_quadrature} provides background on quantization-based cubature formulas which are needed for the following. In Section \ref{seq:generalized_bridges}, we develop a notion of generalized bridge of a continuous Gaussian semimartingale which extends the generalized Brownian bridge introduced by Alili in \cite{AliliGeneralizedBridge}. We prove that under an additional hypothesis (\ref{eq:H_hypothesis}), the generalized bridge of a continuous Gaussian semimartingale remains a Gaussian semimartingale with respect to a bigger filtration and we derive its canonical decomposition. (Let us mention the thorough study of the properties of Gaussian semimartingales available in \cite{JainMonradGaussQuasimart}.) A similar result is stated when conditioning by a Voronoi quantizer. We pay a particular attention to the special case of generalized bridges that we call Karhunen-Loève generalized bridges and which amounts to the conditioning of $X$ by a finite subset of its K-L coordinates. Section \ref{sec:KLgeneralizedBridges} is devoted to the partial functional quantization of continuous Gaussian semimartingales and its application to the partial quantization of solutions of SDEs. We finally give $L^p$ and $a.s.$ convergence results for partially quantized SDEs.

\section{Quantization-based cubature and related inequalities}\label{seq:quantization_based_quadrature}
\par The idea of quantization-based cubature method is to approach the probability distribution of the random variable $X$ by the distribution of a quantizer $Y$ of $X$. As $Y$ is a discrete random variable, we can write $\PP_Y = \sum\limits_{i=1}^N p_i \delta_{y_i}$. If $F : E \to \R$ is a Borel functional, 
\begin{equation}\label{eq:discrete_expectation}
\E[F(Y)] = \sum\limits_{i=1}^N p_i F(y_i).
\end{equation}
Hence, the weighted discrete distribution $(y_i,p_i)_{1 \leq i \leq N}$ of $Y$ allows one to compute the sum (\ref{eq:discrete_expectation}). We review here some error bounds which can be derived when approaching $\E[F(X))]$ by (\ref{eq:discrete_expectation}). See \cite{PagesPrintemsFunctional4} for detailed proofs.
\begin{enumerate}
\item \textit{ If $X \in L^2$, $Y$ a quantizer of $X$ of size $N$ and $F$ is Lipschitz continuous, then 
\begin{equation}\label{eq:error_first_order}
|\E[F(X)] - \E[F(Y)]| \leq [F]_{\textnormal{Lip}} \|X-Y\|_2,
\end{equation}
where $[F]_{\textnormal{Lip}}$ is the Lipschitz constant of $F$. In particular, if $(Y_N)_{N\geq 1}$ is a sequence of quantizers such that $\lim\limits_{N\to \infty} \|X-Y_N\|_2 = 0$, then the distribution $\sum\limits_{i=1}^N p^N_i \delta_{x_i^N}$ of $Y_N$ weakly converges to the distribution $\PP_X$ of $X$ as $N \to \infty$.
}
\vspace{1mm}
\par \noindent \hspace{-10mm} This first error bound is a straightforward consequence of $|F(X)-F(Y)| \leq [F]_{\textnormal{Lip}}|X-Y|$.
\item \textit{If $Y$ is a stationary quantizer of $X$, \textit{i.e.} $Y = \E[X|Y]$, and $F$ is differentiable with an $\alpha$-Hölder differential $DF$ for $\alpha \in (0,1]$, \textit{i.e.} $|DF(u) - DF(v)|_{L(E)} \leq [DF]_\alpha |u-v|^\alpha$, for all $(u,v) \in E^2$ where $|\cdot|_{L(E)}$ is the operator norm on $L(E)$, then
\begin{equation}\label{eq:error_holder_derivative}
|\E[F(X)] - \E[F(Y)]| \leq [DF]_\alpha \|X-Y\|_2^{1+\alpha}.
\end{equation}
In the case where $F$ has a Lipschitz continuous derivative $(\alpha = 1)$, we have. $[DF]_1=[DF]_\textnormal{Lip}$. For example, if $F$ is twice differentiable and $D^2F$ is bounded, then $[DF]_{\textnormal{Lip}} = \|D^2 F\|_\infty$.}
\vspace{1mm}
\par \noindent \hspace{-10mm} This particular inequality comes from the Taylor expansion of $F$ around $X$ and the stationarity of $Y$.
\item \textit{If $F$ is a semi-continuous\footnote{In the infinite-dimensional case, convexity does not imply continuity. In infinite-dimensional Banach spaces, a semi-continuity hypothesis is required for Jensen's inequality. See \cite{ZapalaJensen} for more details.}convex functional and $Y$ is a stationary quantizer of $X$, 
\begin{equation}\label{eq:error_convex}
\E[F(Y)] \leq \E[F(X)].
\end{equation}}
\vspace{-4mm}
\par \noindent \hspace{-10mm} This inequality is a straightforward consequence of the stationarity property and Jensen's inequality.
$$
\E[F(Y)] = \E[F(\E[X|Y])] \leq \E[\E[F(X)|Y]] = \E[F(X)].
$$
\end{enumerate}
\section{Functional quantization and generalized bridges}\label{seq:generalized_bridges}
\subsection{Generalized bridges}\label{seq:generalized_bridge_intro}
\par \noindent Let $(X_t)_{t \in [0,T]}$ be a continuous centered Gaussian semimartingale starting from $0$ on $(\Omega,\mathcal{A},\PP)$ and $\mathcal{F}^X$ its natural filtration. Fernique's theorem ensures that $\int_0^T \E\left[X_t^2\right] dt < +\infty$ (see Janson \cite{JansonGaussianHilbertSpaces}).
\par \noindent We aim here to compute the conditioning with respect to a finite family $\overline{Z}_T:=(Z_T^i)_{i \in I}$ of Gaussian random variables, which are measurable with respect to $\sigma(X_t,t \in[0,T])$. ($I \subset \N$ is a finite subset of $\N^*$.) As Alili in \cite{AliliGeneralizedBridge} we settle on the case where $(Z_T^i)_{i \in I}$ are the terminal values of processes of the form $Z_t^i = \int_0^t f_i(s) dX_s$, $i \in I$, for some given finite set $\overline{f} = (f_i)_{i \in I}$ of $L^2_{loc}([0,T])$ functions. The \emph{generalized bridge} for $(X_t)_{t \in [0,T]}$ corresponding to $\overline{f}$ with end-point $\overline{z} = (z_i)_{i\in I}$ is the process $\left(X_t^{\overline{f},\overline{z}}\right)_{t \in [0,T]}$ that has the distribution
\begin{equation}\label{eq:generalized_bridge}
X^{\overline{f},\overline{z}} \simdist \mathcal{L}\left(X \middle| Z_T^i = z_i, \ i\in I \right). 
\end{equation}
\par \noindent For example, in the case where $X$ is a standard Brownian motion with $|I| = 1$, $\overline{f} = \{ f \}$ and $f\equiv 1$, this is the Brownian bridge on $[0,T]$. If $X$ is an Ornstein-Uhlenbeck process this is an Ornstein-Uhlenbeck bridge.
\vspace{2mm}
\par \noindent Let $H$ be the Gaussian Hilbert space spanned by $(X_s)_{s\in [0,T]}$ and $H_{\overline{Z}_T}$ the closed subspace of $H$ spanned by $(Z^i_T)_{i \in I}$. We denote by $H_{\overline{Z}_T}^\perp$ its orthogonal complement in $H$. Any Gaussian random variable $G$ of $H$ can be orthogonally decomposed into $G = \Proj_{\overline{Z}_T}(G) {\overset{\independent }{+}} \Proj_{\overline{Z}_T}^\perp(G),$ where $\Proj_{\overline{Z}_T}$ and $\Proj_{\overline{Z}_T}^\perp$ are the orthogonal projections on $H_{\overline{Z}_T}$ and $H_{\overline{Z}_T}^\perp$. ($\Proj_{\overline{Z}_T}^\perp = Id_H - \Proj_{\overline{Z}_T}$). With these notations, $\E\left[G\middle|(Z_T^i)_{i \in I}\right] = \Proj_{\overline{Z}_T}(G)$. 
\par Other definitions of generalized bridges exist in the literature, see \textit{e.g.} \cite{YorMansuy}. 
\subsection{The case of the Karhunen-Loève basis}\label{sec:KL_generalized_bridge}
\par \noindent As $X$ is a continuous Gaussian process, it has a continuous covariance function (see \cite[VIII.3]{JansonGaussianHilbertSpaces}). We denote by $(e_i^X,\lambda_i^X)_{i \geq 1}$ its Karhunen-Loève eigensystem. Thus, if we define function $f^X_i$ as the antiderivative of $-e_i^X$ that vanishes at $t=T$, \textit{i.e.} $f^X_i(t) = \int_t^T e_i^X(s)ds$, an integration by parts yields
\begin{equation}\label{eq:primitive_KL}
\int_0^T X_s e_i^X(s) ds = \int_0^T f^X_i(s) dX_s.
\end{equation}
\par \noindent In other words, with the notations of Section \ref{seq:generalized_bridge_intro}, we have $Z_T^i = \int_0^T X_s e_i^X(s) ds =: Y_i$,  the $i$th Karhunen-Loève coordinate of $X$.
\par \noindent For some finite subset $I \subset \N^*$, we denote by $X^{I,\overline{y}}$ and call \emph{K-L generalized bridge} the generalized bridge associated with functions $(f_i^X)_{i \in I}$ and with end-point $\overline{y} = (y_i)_{i\in I}$. This process has the distribution $\mathcal{L}(X|Y_i = y_i, i \in I)$.
\par \noindent In this case, the Karhunen-Loève expansion gives the decomposition
\begin{equation}\label{eq:conditioning_KL}
X = \underbrace{\sum\limits_{i \in I} Y_i e_i^X}_{=\Proj_{\overline{Z}_T}(X)} {\overset{\independent }{+}} \underbrace{\sum\limits_{i \in \N^* \backslash I} \sqrt{\lambda_i^X} \xi_i e_i^X}_{=\Proj_{\overline{Z}_T}^\perp(X)},
\end{equation}
where $(\xi_i)_{i \in \N^* \backslash I}$ are independent standard Gaussian random variables. This gives us the projections $\Proj_{\overline{Z}_T}$ and $\Proj_{\overline{Z}_T}^\perp$ defined in Section \ref{seq:generalized_bridge_intro}. It follows from (\ref{eq:conditioning_KL}) that a K-L generalized bridge is centered on $\E\left[X \middle|Y_i = y_i, i \in I \right]$ and has the covariance function 
\begin{equation}\label{eq:conditional_covariance_function}
\Gamma^{X|Y}(s,t) = \cov(X_s,X_t) - \sum\limits_{i \in I} \lambda_i^X e_i^X(s) e_i^X(t).
\end{equation}
\par \noindent We have $\int_0^T \Gamma^{X|Y}(t,t) dt = \sum\limits_{i \in \N^* \backslash I} \lambda_i^X$.
\par \noindent Moreover, thanks to Decomposition (\ref{eq:conditioning_KL}), if $X^{I,\overline{y}}$ is a K-L generalized bridge associated with $X$ with terminal values $\overline{y}=(y_i)_{i \in I}$, it has the same probability distribution as the process 
$$
\sum\limits_{i\in I} y_i e_i^X(t) + X_t - \sum\limits_{i \in I} \left( \int_0^T X_s e_i^X(s) ds \right) e_i^X(t).
$$
\par \noindent This process is then the sum of a semimartingale and a non-adapted finite-variation process. 
\par \noindent Let us stress the fact that the second term in the left-hand side of (\ref{eq:conditioning_KL}) is the corresponding K-L generalized bridge with end-point $0$, \textit{i.e.} $\Proj_{\overline{Z}_T}^\perp = X^{I,\overline{0}}$.
\par \noindent In \cite{CorlayPagesStratification}, an algorithm is proposed to exactly simulate marginals of a K-L generalized bridge with a linear additional cost to a prior simulation of $(X_{t_0},\cdots,X_{t_n})$, for some subdivision $0 = t_0 \leq t_1 \leq \cdots t_n = T$ of $[0,T]$. This was used for variance reduction issues. Note that the algorithm is easily extended to the case of (non-K-L) generalized bridges. 
\subsection{Generalized bridges as semimartingales}\label{sec:semimartingale_generalized}
\par \noindent For a random variable $L$, we denote by $\PP\left[\cdot \middle|L \right]$ the conditional probability knowing $L$. We keep the notations and assumptions of previous sections. ($X$ is a continuous Gaussian semimartingale starting from $0$.) We consider a finite set $I \subset \{1, 2, \cdots \}$ and $(f_i)_{i\in I}$ a set of bounded measurable functions. Let $X^{\overline{f},\overline{z}}$ be the generalized bridge associated with $X$ with end-point $\overline{z} = (z_i)_{i \in I}$. For $i \in I$, $Z_t^i = \int_0^t f_i(s) dX_s$ and $\overline{Z}_t = (Z_t^i)_{i \in I}$. 
\vspace{3mm}
\par \noindent Jirina's theorem ensures the existence of a transition kernel
$$
\nu_{\left.\overline{Z}_T \middle| \left((X_t)_{t \in [0,s]}\right) \right.} : \mathcal{B}(\R^I) \times C^0\left([0,s],\R \right) \to \R_+,
$$
corresponding to the conditional distribution $\mathcal{L}\left(\overline{Z}_t \middle| \left((X_t)_{t \in [0,s]}\right)\right)$. 
\par \noindent We now make the additional assumption (\ref{eq:H_hypothesis}) that, for every $s \in [0,T)$ and for every $(x_u)_{u \in [0,s]} \in C^0\left([0,s],\R \right)$, the probability measure $\nu_{\left.\overline{Z}_T \middle| \left((X_t)_{t \in [0,s]}\right) \right.}\left(d\overline{z},(x_u)_{u \in [0,s]}\right)$ is absolutely continuous with respect to the Lebesgue measure. We denote by $\Pi_{(x_u)_{u \in [0,s]},T}$ its density. The covariance matrix of this Gaussian distribution on $\R^I$ writes
$$
Q(s,T) := \E\left[\left(\overline{Z}_T - \E\left[\overline{Z}_T\middle|(X_u)_{u \in [0,s]}\right] \right)\left(\overline{Z}_T - \E\left[\overline{Z}_T \middle|(X_u)_{u \in [0,s]}\right] \right)^* \middle| (X_u)_{u \in [0,s]}\right].
$$
\par \noindent If $X$ is a martingale, we have $Q(s,T) = \left( \left( \int_s^T f_i(u)f_j(u) d \langle X \rangle_u \right)\right)_{(i,j)\in I^2}$. We recall that a continuous centered semimartingale $X$ is Gaussian if and only if $\langle X \rangle$ is deterministic (see \textit{e.g.} \cite{RevuzYor}). Hence, this additional hypothesis is equivalent to assume that 
\begin{equation}\tag{$\mathcal{H}$}\label{eq:H_hypothesis}
Q(s,T) \quad \textnormal{is invertible for every } s\in [0,T).
\end{equation}
\par \noindent The following theorem follows from the same approach as the homologous result in the article by Alili \cite{AliliGeneralizedBridge} for the Brownian case. It is extended to the case of a continuous centered Gaussian semimartingale starting from $0$. 
\begin{theo}\label{thm:alili}
\vspace{2mm}
\par \noindent Under the (\ref{eq:H_hypothesis}) hypothesis, for any $s \in [0,T)$, and for $\PP_{\overline{Z}_T}$-almost every $\overline{z}\in \R^I$, $\PP\left[\cdot \middle|\overline{Z}_T = \overline{z} \right]$ is equivalent to $\PP$ on $\mathcal{F}^X_s$ and its Radon-Nikodym density is given by
$$
\frac{d\PP\left[\cdot \middle|\overline{Z}_T = \overline{z} \right]}{d\PP}_{| \mathcal{F}^X_s} = \frac{\Pi_{(X_u)_{u\in [0,s]},T}(\overline{z})}{\Pi_{0,T}(\overline{z})}.
$$
\end{theo}
\par \noindent \textbf{Proof:} Consider $F$ a real bounded $\mathcal{F}^X_s$-measurable random variable and $\phi:\R^I \to \R$ a bounded Borel function. 
\begin{itemize}
\item On the one hand, preconditioning by $\overline{Z}_T$ yields
\begin{equation}\label{eq:bayes_bridge_1}
\E\left[F \phi(\overline{Z}_T)\right] = \E\left[ \E\left[F \middle| \overline{Z}_T \right] \phi(\overline{Z}_T) \right] = \int\limits_{\R^I} \phi(\overline{z}) \E\left[F\middle|\overline{Z}_T = \overline{z}\right] \Pi_{0,T}(\overline{z})d\overline{z}. 
\end{equation}
\item On the other hand, as $F$ is measurable with respect to $\mathcal{F}^X_s$, preconditioning with respect to $\mathcal{F}^X_s$ yields
$$
\E\left[F \phi(\overline{Z}_T)\right] = \E\left[F \E \left[\phi(\overline{Z}_T)\middle| \mathcal{F}^X_s \right]\right] = \E\left[ F \int_{\R^I} \phi(\overline{z}) \Pi_{(X_t)_{t\in [0,s]},T}(\overline{z}) d\overline{z} \right].
$$
\par \noindent Now, thanks to Fubini's theorem
\begin{equation}\label{eq:bayes_bridge_2}
\E\left[F \phi(\overline{Z}_T)\right] = \int\limits_{\R^I} \phi(\overline{z}) \E\left[F \Pi_{(X_t)_{t\in [0,s]},T}(\overline{z})\right] d\overline{z}.
\end{equation}
\end{itemize}
Identifying Equations (\ref{eq:bayes_bridge_1}) and (\ref{eq:bayes_bridge_2}), we see that for $\PP_{\overline{Z}_T}$-almost surely $\overline{z}\in \R^I$ and for every real bounded $\mathcal{F}^X_s$-measurable random variable $F$, 
\begin{equation}\label{eq:radon_nikodym_characterization}
\E\left[F \middle| \overline{Z}_T = \overline{z}\right] = \E\left[F \frac{\Pi_{(X_t)_{t\in [0,s]},T}(\overline{z})}{\Pi_{0,T}(\overline{z})}\right].
\end{equation}
\par \noindent Equation (\ref{eq:radon_nikodym_characterization}) characterizes the Radon-Nikodym derivative of $\PP\left[\cdot \middle|\overline{Z}_T = \overline{z}\right]$ with respect to $\PP$ on $\mathcal{F}^X_s$. \myqed
\vspace{5mm}
\par \noindent We now can use classical filtration enlargement techniques \cite{JacodEnlargement,JeulinEnlargement,YorEnlargement}. 
\begin{prop}[Generalized bridges as semimartingales]\label{prop:bridge_as_semimartingale}
\par \noindent Let us define the filtration $\mathcal{G}^{X,\overline{f}}$ by $\mathcal{G}^{X,\overline{f}}_t := \sigma\left(\overline{Z}_T, \mathcal{F}^X_t\right)$, the enlargement of the filtration $\mathcal{F}^X$ corresponding to the above conditioning. We consider the stochastic process $D_s^{\overline{z}} := \frac{d \PP\left[\cdot\middle|\overline{Z}_T = \overline{z} \right]}{d\PP}_{|\mathcal{F}_s^X} = \frac{\Pi_{(X_t)_{t\in [0,s]},T}(\overline{z})}{\Pi_{0,T}(\overline{z})}$ for $s \in [0,T)$. 
\par Under the (\ref{eq:H_hypothesis}) hypothesis, and the assumption that $D^{\overline{z}}$ is continuous, $X$ is a continuous $\mathcal{G}^{X,\overline{f}}$-semimartingale on $[0,T)$. 
\end{prop}
\par \noindent \textbf{Proof:} $D^{\overline{z}}$ is a strictly positive martingale on $[0,T)$ which is uniformly integrable on every interval $[0,t] \subset [0,T)$. Hence, as we assumed that it is continuous, we can write $D^{\overline{z}}$ as an exponential martingale $D_s^{\overline{z}} = \exp\left(L_s^{\overline{z}} - \frac{1}{2} \left\langle L^{\overline{z}} \right\rangle_s \right)$ with $L_t^{\overline{z}} = \int_0^t \left(D_s^{\overline{z}}\right)^{-1} dD_s^{\overline{z}}$ (as $D_0^{\overline{z}}=1$).
\par \noindent Now, as $X$ is a continuous $(\mathcal{F}^X,\PP)$-semimartingale, we write $X = V+M$ its canonical decomposition (under the filtration $\mathcal{F}^X$). 
\begin{itemize}
\item Thanks to Girsanov theorem, $\widetilde{M}^{\overline{z}}:= M - \left\langle M,L^{\overline{z}}\right\rangle$ is a $\left(\mathcal{F}^X,\PP\left[\cdot\middle|\overline{Z}_T = \overline{z}\right]\right)$-martingale. 
\begin{itemize}
\item A consequence is that it is a $\left(\mathcal{G}^{X,\overline{f}},\PP\left[\cdot \middle|\overline{Z}_T= \overline{z}\right]\right)$-martingale.
\item And thus $\widetilde{M}^{\overline{Z}_T}$ is a $\left(\mathcal{G}^{X,\overline{f}},\PP\right)$-martingale.
\par \noindent For more preciseness on this, we refer to \cite[Theorem $3$]{AnkirchnerDereichImkeller} where the proof is based on the notion of decoupling measure.
\end{itemize}
\item Moreover, conditionally to $\overline{Z}_T$, $V$ is still a finite-variation process $V$, and is adapted to $\mathcal{G}^{X,\overline{f}}$. \myqed
\end{itemize}
\begin{remark}[Continuous modification]
\par \noindent In Proposition \ref{prop:bridge_as_semimartingale}, if one only assumes that $D^{\overline{z}}$ has a continuous modification $\mathcal{D}^{\overline{z}}$, then with each one of its continuous modifications is associated a continuous $\mathcal{G}^{X,\overline{f}}$-semimartingale on $[0,T)$, and all these semimartingales are modifications of each other. 
\end{remark}
\begin{prop}[Continuity of $D^{\overline{z}}$]
\par If $\mathcal{F}^X$ is a standard Brownian filtration, then $D^{\overline{z}}$ has a continuous modification. 
\end{prop}
\par \noindent \textbf{Proof:} Consider $s \in [0,T)$. Under the (\ref{eq:H_hypothesis}) hypothesis, the density $\Pi_{(X_u)_{u\in [0,s]},T}$ writes
{\footnotesize
\begin{equation}\label{eq:bi_density}
\Pi_{(X_u)_{u\in [0,s]},T}(\overline{z}) = \left( 2 \pi \det Q(s,T) \right)^{-\frac{|I|}{2}} \exp \left(\left(\overline{z}- \E\left[ \overline{Z}_T \middle| (X_u)_{u \in [0,s]} \right] \right) Q(s,T)^{-1} \left(\overline{z}-\E\left[ \overline{Z}_T \middle| (X_u)_{u \in [0,s]} \right] \right)^* \right).
\end{equation}
}
\par \noindent Let us define the stochastic process $\overline{H}$ by $\overline{H}_s := \E\left[ \overline{Z}_T \middle| (X_u)_{u \in [0,s]} \right]$. The so-defined process $\overline{H}$ is a $\mathcal{F}^X$ local martingale. Thanks the Brownian representation theorem, $\overline{H}$ has a Brownian representation and has a continuous modification. The continuity of $s \mapsto \det Q(s,T)$ and $s \mapsto Q(s,T)^{-1}$ follows from the definition of $Q(s,T)$ and the continuity of $\overline{H}$ (up to a modification). Hence, $D^{\overline{z}}$ has a continuous modification. \myqed
\begin{remark}
\begin{itemize}
\item The measurability assumption with respect to a Brownian filtration is satisfied in the cases of Brownian bridge and Ornstein-Uhlenbeck processes.
\item This hypothesis is not necessary so long as the continuity of the martingale $\overline{H}_s = \E\left[ \overline{Z}_T \middle| (X_u)_{u \in [0,s]} \right]$ can be proved by any means.
\end{itemize}
\end{remark}
\subsubsection{On the canonical decomposition}
\par \noindent With the same notations, and under the (\ref{eq:H_hypothesis}) hypothesis, we can tackle the canonical decomposition of $X^{\overline{f},\overline{z}}$. We have
$$
L_t^{\overline{z}} = \int_0^t \frac{d\Pi_{(X_u)_{u\in [0,s]},T}(\overline{z})}{\Pi_{(X_u)_{u\in [0,s]},T}(\overline{z})},
$$
and 
\begin{multline*}
\ln \left(\Pi_{(X_u)_{u\in [0,s]},T}(\overline{z})\right) = -\frac{|I|}{2}\ln \left(2 \pi \det Q(s,T) \right) \\
- \frac{1}{2} \left(\overline{z} - \E \left[\overline{Z}_T \middle|(X_u)_{u \in[0,s]} \right] \right) Q(s,T)^{-1} \left(\overline{z} - \E \left[\overline{Z}_T \middle|(X_u)_{u \in[0,s]} \right] \right)^*.
\end{multline*}
\par \noindent Using that for a positive continuous semimartingale $S$, $d\ln S = \frac{dS}{S}-\frac{1}{2}d \left\langle\frac{1}{S} \cdot S\right\rangle$, we obtain
$$
\begin{array}{lll}
\frac{d\Pi_{(X_u)_{u\in [0,s]},T}(\overline{z})}{\Pi_{(X_u)_{u\in [0,s]},T}(\overline{z})} &= d \ln \left(\Pi_{(X_u)_{u\in [0,s]},T}(\overline{z}) \right) + \left(\begin{array}{cc} \textnormal{finite-variation}\\ \textnormal{process} \end{array}\right) \\
				&=-\frac{1}{2} d \left(\left(\overline{z}-\E \left[\overline{Z}_T \middle|(X_u)_{u \in[0,s]} \right] \right) Q(s,T)^{-1} \left(\overline{z}-\E \left[\overline{Z}_T \middle|(X_u)_{u \in[0,s]} \right] \right)^*\right) + \left(\textnormal{f.-v. p.}\right) \\
				& = \left(d \E \left[\overline{Z}_T \middle|(X_u)_{u \in[0,s]} \right] \right) Q(s,T)^{-1} \left(\overline{z}-\E \left[\overline{Z}_T \middle|(X_u)_{u \in[0,s]} \right] \right)^*+\left(\textnormal{f.-v. p.}\right).
\end{array}
$$
Hence,
$$
d \left\langle X, L^{\overline{z}} \right\rangle_s = d \left\langle X, \E \left[\overline{Z}_T \middle|(X_u)_{u \in[0,\cdot]} \right] \right\rangle_s Q(s,T)^{-1} \left(\overline{z}-\E \left[\overline{Z}_T \middle|(X_u)_{u \in[0,s]} \right] \right)^*.
$$
This expression can be further simplified in the two following cases:
\begin{itemize}
\item \emph{In the case where $X$ is a martingale}, owing to the definition of $Z_j$, we have $\forall j \in I, \ \E \left[Z^j_T \middle|(X_u)_{u \in[0,s]} \right] = \int_0^s f_j(u) dX_u$ so that
\begin{equation}\label{eq:final_canonical_decomposition}
\begin{array}{lll}
d \left\langle X, L^{\overline{z}} \right\rangle_s 	& = \left(\overline{f}(s) Q(s,T)^{-1} \left(\overline{z}- \E \left[\overline{Z}_T \middle|(X_u)_{u \in[0,s]} \right] \right)^* \right) d\langle X \rangle_s\\
							& = \sum\limits_{i \in I} f_i(s) \sum\limits_{j\in I} \left(Q(s,T)^{-1}\right)_{ij} \left(z_j - \E \left[Z^j_T \middle|(X_u)_{u \in[0,s]} \right] \right) d\langle X \rangle_s.
\end{array}
\end{equation}
\par \noindent As a consequence, $M - \int_0^\cdot \sum\limits_{i \in I} f_i(s) \sum\limits_{j\in I} \left(Q(s,T)^{-1}\right)_{ij} \left(z_j - \E \left[Z^j_T \middle|(X_u)_{u \in[0,s]} \right] \right) d\langle X \rangle_s$ is a $\left(\mathcal{G}^{X,\overline{f}},\PP\left[\cdot \middle| \overline{Z}_T = \overline{z}\right] \right)$-martingale. We have recovered Alili's result on the generalized Brownian bridge \cite{AliliGeneralizedBridge}. 
\item \emph{In the case where the Gaussian semimartingale $X$ is a Markov process}, for every $j \in I$ there exists $g_j \in L^2 ([0,T])$ such that $\E \left[Z^j_T \middle|(X_u)_{u \in[0,s]} \right] = \int_0^s f_j(u) dX_u + g_j(s) X_s$. Indeed, 
$$
\E \left[Z^j_T \middle|(X_u)_{u \in[0,s]} \right] = \int_0^s f_j(u) dX_u + \underbrace{\E \left[\int_s^T f_j(u) dX_u \middle|(X_u)_{u \in[0,s]} \right]}_{=:g_j(s) X_s}.
$$
Hence, if one assumes that $(g_j)_{j \in I}$ are finite-variation functions (which is the case when $X$ is an Ornstein-Uhlenbeck process or a Brownian bridge), we have $d \left\langle X, \E \left[\overline{Z}_T \middle|(X_u)_{u \in[0,\cdot]} \right]\right\rangle_s = \left( \overline{f}(s) + \overline{g}(s) \right) d\langle X \rangle_s$, and thus
$$
\begin{array}{lll}
d \left\langle X, L^{\overline{z}} \right\rangle_s 	& = \left(\left(\overline{f}(s)+\overline{g}(s)\right) Q(s,T)^{-1} \left(\overline{z}- \E \left[Z^j_T \middle|(X_u)_{u \in[0,s]} \right] \right)^* \right) d\langle X \rangle_s\\
							& = \sum\limits_{i \in I} \left(f_i(s) + g_i(s) \right) \sum\limits_{j\in I} \left(Q(s,T)^{-1}\right)_{ij} \left(z_j - \E \left[Z^j_T \middle|(X_u)_{u \in[0,s]} \right] \right) d\langle X \rangle_s.
\end{array}
$$
\end{itemize}
\begin{example}[Standard Brownian bridge]
In the case where $X = W$ is a standard Brownian motion with $|I| = 1$, $\overline{f} = \{ f \}$ and $f\equiv 1$, $Z_t = W_t$ and $W^{\overline{f}, \overline{z}}$ is a standard Brownian bridge. We have $Q(s,T)^{-1} = \frac{1}{T-s}$ and
$$
d \langle X, L^z \rangle_t = \frac{1}{T-t} \left( z - \E\left[W_T \middle| (W_u)_{u \in [0,t]}\right]\right)dt = \frac{z-W_t}{T-t} dt. 
$$
Thus, 
$$
dW_t = \frac{z-W_t}{T-t}dt \quad + \!\!\!\! \underbrace{\frac{W_t}{T-t}dt + dW_t.}_{\left( \mathcal{G}^{X,\overline{f}}, \PP\left[\cdot \middle| W_T = z \right]\right)\textnormal{-martingale}}
$$
The martingale part happens to be a $\left( \mathcal{G}^{X,\overline{f}}, \PP\left[\cdot \middle| W_T = z \right] \right)$-standard Brownian motion, thanks to Lévy's characterization of the Brownian motion. Thus we have retrieved the classical SDE of the Brownian bridge. 
\end{example}
\subsubsection{Generalized bridges and functional stratification}
\par \noindent With the same notations, we set $\widehat{Z}^\Gamma = \Proj_\Gamma(\overline{Z}_T) = \sum\limits_{i=1}^N \gamma_i \1_{C_i}(\overline{Z}_T)$ a stationary quantizer of $\overline{Z}_T$ (where $\Gamma = \{\gamma_1,\cdots,\gamma_N\}$ and $C = \{C_1,\cdots,C_N\}$ are respectively the associated knots and Voronoi partition). 
\begin{prop}[Stratification]\label{prop:stratif_equivalent}
\par Under the (\ref{eq:H_hypothesis}) hypothesis, for any $s \in [0,T)$, for any $k \in \left\{ 1, \cdots, N \right\}$, $\PP\left[\widehat{Z}^\Gamma = \gamma_k\right]>0$ and the conditional probability $\PP\left[\cdot \middle| \widehat{Z}^\Gamma = \gamma_k \right]$ is equivalent to $\PP$ on $\mathcal{F}^X_s$. 
\end{prop}
\par \noindent \textbf{Proof:} Obviously, if $A\in \mathcal{F}^X_s$ is such that $\PP[A] = 0$, we have $\PP\left[A \middle| \widehat{Z}^\Gamma = \gamma_k \right] = 0$. Conversely, $B \in \mathcal{F}^X_s$ satisfies $\PP\left[B \middle| \widehat{Z}^\Gamma = \gamma_k\right] = 0$, then pre-conditioning by $\overline{Z}_T$, we get $\E \left[ \E \left[ \1_B \middle| \overline{Z}_T \right] \middle| \widehat{Z}^\Gamma = \gamma_k \right]= 0.$ Thus, $\int_{\overline{z}\in C_k} \PP\left[B \middle| \overline{Z}_T = \overline{z} \right] d \PP_{\overline{Z}_T}(\overline{z})= 0$. Hence $\PP\left[B \middle| \overline{Z}_T = \overline{z} \right] =0$ for $\PP_{\overline{Z}_T}$-almost every $\overline{z} \in C_k$. 
\par \noindent Since $\PP_{\overline{Z}_T}(C_k)>0$, there exists at least one element $\overline{z} \in C_k$ such that $\PP\left[B \middle| \overline{Z}_T = \overline{z} \right] =0$. Now thanks to Theorem \ref{thm:alili}, $\PP[B]=0$. \myqed
\begin{prop}[Stratification]\label{prop:stratif_semimartingale}
\par \noindent Let us define the filtration $\mathcal{G}^{X,\Gamma}$ by $\mathcal{G}^{X,\Gamma}_t := \sigma\left(\mathcal{F}^X_t,\widehat{Z}^{\Gamma}\right)$, the enlargement of $\mathcal{F}^X$ corresponding to the conditioning with respect to $\widehat{Z}^{\Gamma}$. For $k \in \left\{1,\cdots,N\right\}$, we consider the stochastic process $D^{\gamma_k}_s := \frac{d\PP\left[\cdot \middle| \widehat{Z}^\Gamma = \gamma_k \right]}{d\PP}_{|\mathcal{F}_s^X}$ for $s \in [0,T)$. 
\par Under the (\ref{eq:H_hypothesis}) hypothesis, and the assumption that $D^{\gamma_k}$ is continuous, the conditional distribution $\mathcal{L}\left(X \middle|\widehat{Z}^{\Gamma} \right)$ of $X$ knowing in which Voronoi cell $\overline{Z}_T$ falls, is the probability distribution of a $\mathcal{G}^{X,\Gamma}$-semimartingale on $[0,T)$.
\end{prop}
\par \noindent \textbf{Proof:} Using that $\PP\left[\cdot \middle| \widehat{Z}^\Gamma = \gamma_k \right]$ is equivalent to $\PP$ on $\mathcal{F}^X_s$, thanks to Proposition \ref{prop:stratif_equivalent}, we can \textit{mutatis mutandis} use the same arguments as for Proposition \ref{prop:bridge_as_semimartingale}, $\PP\left[ \cdot \middle| \overline{Z}_T = \overline{z} \right]$ being replaced by $\PP\left[\cdot\middle|\widehat{Z}^\Gamma = \gamma_k \right]$.
\par $D^{\gamma_k}$ is a strictly positive martingale on $[0,T)$ uniformly integrable on every $[0,t] \subset [0,T)$. Hence, as $D^{\gamma_k}$ is continuous by hypothesis, it is an exponential martingale $D^{\gamma_k}_s = \exp\left(L^{\gamma_k}_s - \frac{1}{2} \left\langle L^{\gamma_k} \right\rangle_s\right)$, with $L^{\gamma_k}_t = \int_0^t \left(D^{\gamma_k}_s\right)^{-1} dD^{\gamma_k}_s$ (as $D^{\gamma_k}_0 = 1$). Now, as $X$ is a continuous $(\mathcal{F}^X,\PP)$-semimartingale, we write $X = V+M$ its canonical decomposition (under the filtration $\mathcal{F}^X$). 
\begin{itemize}
\item Thanks to Girsanov theorem, $\widetilde{M}^{\gamma_k} := M - \left\langle M,L^{\gamma_k} \right\rangle$ is a $\left(\mathcal{F}^X,\PP\left[\cdot \middle|\widehat{Z}^\Gamma = \gamma_k \right]\right)$-martingale. 
As a consequence, it is a $\left(\mathcal{G}^{X,\Gamma},\PP\left[\cdot\middle|\widehat{Z}^\Gamma = \gamma_k \right]\right)$-martingale and thus $\widetilde{M}^{\widehat{Z}^\Gamma}$ is a $\left(\mathcal{G}^{X,\Gamma},\PP\right)$-martingale. 
\item Moreover, conditionally to $\widehat{Z}^\Gamma$, $V$ is still a finite-variation process $V$, and is adapted to $\mathcal{G}^{X,\Gamma}$. \myqed
\end{itemize}
\begin{prop}[Continuity of $D^{\gamma_k}$]
\par If $\mathcal{F}^X$ is a Brownian filtration, then $D^{\gamma_k}$ has a continuous modification. 
\end{prop}
\par \noindent \textbf{Proof:} By definition, $D^{\gamma_k}$ is a $\mathcal{F}^X$-local martingale on $[0,T]$. The conclusion is a straightforward consequence of the Brownian representation theorem.\myqed
\vspace{2mm}
\par \noindent Considering the partition of $L^2([0,T])$ corresponding to the Voronoi cells of a functional quantizer of $X$, the last two propositions show that the conditional distribution of the $X$ in each Voronoi cell (strata) is a Gaussian semimartingale with respect to its own filtration. This allows us to define the corresponding functional stratification of the solutions of stochastic differential equations driven by $X$. 
\vspace{2mm}
\par \noindent In \cite{CorlayPagesStratification}, an algorithm is proposed to simulate the conditional distribution of the marginals $(X_{t_0},\cdots,X_{t_n})$ of $X$ for a given subdivision $0 = t_0 < t_1 < \cdots < t_n = T$ of $[0,T]$ conditionally to a given Voronoi cell (strata) of a functional quantization of $X$. The simulation complexity has an additional linear complexity to an unconditioned simulation of $(X_{t_0},\cdots,X_{t_n})$. We refer to \cite{CorlayPagesStratification} for more details. 
\vspace{2mm}
\par \noindent To deal with the solution of a SDE, it was proposed in \cite{CorlayPagesStratification} to simply plug these marginals in the Euler scheme of the SDE. Proposition \ref{prop:stratif_semimartingale} now shows that this amounts to simulate the Euler scheme of the SDE driven by the corresponding (non-Gaussian) semimartingale. 
\subsection{About the (\ref{eq:H_hypothesis}) hypothesis}
\subsubsection{The martingale case}
\par \noindent In the case where $X$ is a continuous Gaussian martingale, the matrix $Q(s,t)$ defined in Section \ref{sec:semimartingale_generalized} writes $Q(s,t) = \left(\left( \int_s^t f_i(u) f_j(u) d \langle X \rangle_u \right)\right)_{(i,j) \in I^2}$.
\par \noindent For $1 \leq s < t \leq T$, the map $\left( \cdot \middle| \cdot \right) : (f,g) \mapsto \int_s^t f(u) g(u) d\langle X \rangle_u$ defines a scalar product on $L^2([s,t],d\langle X\rangle)$. Hence $Q(s,t)$ is the Gram matrix of the vectors of $L^2([s,t],d\langle X\rangle)$ defined by the restrictions to $[s,t]$ of the functions $(f_i)_{i \in I}$. Thus, it is invertible if and only if these restrictions form a linearly independent family of $L^2([s,t],d\langle X\rangle)$. (Another consequence, is that if $Q(s,t)$ is invertible for some $0 \leq s <t \leq T$, then for every $(u,v)$ such that $[s,t] \subset [u,v]$, $Q(u,v)$ is invertible).
\vspace{2mm}
\par For instance, if $X$ is a standard Brownian motion on $[0,T]$, the functions $(f_i^X)_{i \in I}$ (associated with the Karhunen-Loève decomposition) are trigonometric functions with strictly different frequencies. Hence, they form a linearly independent family of continuous functions on every nonempty interval $[s,T) \subset [0,T)$. Moreover, the measure $d\langle X \rangle$ is proportional to the Lebesgue measure on $[0,T]$ and thus $Q(s,T)$ is invertible for any $s \in [0,T)$. \emph{Hence, the (\ref{eq:H_hypothesis}) hypothesis is fulfilled in the case of K-L generalized bridges of the standard Brownian motion}. 
\subsubsection{Standard Brownian bridge and Ornstein-Uhlenbeck processes}
\par \noindent Brownian bridge and the Ornstein-Uhlenbeck process are not martingales. Hence, this criterion is not sufficient and the invertibility of matrix $Q(s,T)$ has to be proved by other means. 
\par \noindent Following from the definitions of $Q(s,T)$ and $\overline{Z}_T$, in the case of the K-L generalized bridge
\begin{equation}\label{eq:about_H_hyp}
\begin{array}{lll}
Q(s,T)_{ij}	&= \E\Big[\left(\int_s^T f^X_i(u) dX_u - \E\left[\int_s^T f^X_i(u) dX_u \middle|(X_u)_{u \in [0,s]}\right] \right)\\
			&\quad \quad \quad \times \left(\int_s^T f^X_j(u) dX_u - \E\left[\int_s^T f^X_j(u) dX_u \middle|(X_u)_{u \in [0,s]}\right] \right)^* \Big| (X_u)_{u \in [0,s]}\Big]\\
			&= \cov\left(\int_s^T f^X_i(u) dX^{(s)}_u, \int_s^T f^X_j(u) dX^{(s)}_u \right),
\end{array}
\end{equation}
where $\left(X^{(s)}_u\right)_{u \in [s,T]}$ has the conditional distribution of $X$ knowing $(X_u)_{u \in [0,s]}$.
\begin{itemize}
\item When $X$ is a standard Brownian bridge on $[0,T]$, $X_u^{(s)}$ is a Brownian bridge on $[s,T]$, starting from $X_s$ and arriving at $0$.
\par {\small It is the sum of an affine function and a standard centered Brownian bridge on $[s,T]$.}
\item When $X$ is a centered Ornstein-Uhlenbeck process, $X_u^{(s)}$ is an Ornstein-Uhlenbeck process on $[s,T]$ starting from $X_s$, with the same mean reversion parameter as $X$. 
\par {\small It is also the sum of a deterministic function and an Ornstein-Uhlenbeck process starting from $0$.}
\end{itemize}
\par \noindent As a consequence, in these two cases, the quantity $\cov\left(\int_s^T f^X_i(u) dX^{(s)}_u, \int_s^T f^X_j(u) dX^{(s)}_u \right)$ can be computed by plugging either a centered Brownian bridge on $[s,T]$ or an Ornstein-Uhlenbeck starting from $0$ instead of $X^{(s)}$ in Equation (\ref{eq:about_H_hyp}). This means that $Q(s,T)$ is the Gram matrix of the random variables $\left(\int_s^T f^X_i(u) dG_u\right)_{i \in I}$, where the centered Gaussian process $(G_u)_{u \in [s,T]}$ is either a standard Brownian bridge on $[s,T]$ or an Ornstein-Uhlenbeck process starting from $0$ at $s$. Thus it is singular if and only if there exists $(\alpha_i)_{i\in I} \neq 0$ in $\R^I$ such that
\begin{equation}\label{eq:zero_condition}
\int_s^T \underbrace{\left( \sum\limits_{i \in I} \alpha_i f^X_i(u) \right)}_{=: g(u)} dG_u = 0 \quad a.s..
\end{equation}
\par \textbf{The case of Brownian bridge}
\par \noindent In the case where $X$ is the standard Brownian bridge on $[0,T]$, functions $(f^X_i)_{i \in I}$ are $C^\infty$ functions and $G$ is a standard Brownian bridge on $[s,T]$. An integration by parts gives $\int_s^T G_s g'(s) ds = 0 \quad a.s.$ and thus $g' \equiv 0$ on $(s,T)$ and thus $g$ is constant on $[s,T]$. The functions $(f_i^X)_{i \in I}$ form a linearly independent set of functions and, as they are trigonometric functions with different frequencies, they clearly don't span constant functions, so that Equation (\ref{eq:zero_condition}) yields $\alpha_1 = \cdots = \alpha_n = 0$. \emph{Hence the (\ref{eq:H_hypothesis}) hypothesis is fulfilled in the case of K-L generalized bridges of the standard Brownian bridge.}
\vspace{2mm}
\par \textbf{The case of Ornstein-Uhlenbeck processes}
\par \noindent In the case where $X$ is an Ornstein-Uhlenbeck process on $[0,T]$, $G$ is an Ornstein-Uhlenbeck process on $[s,T]$ starting from $0$. The injectivity property of the Wiener integral related to the Ornstein-Uhlenbeck process stated in Proposition \ref{prop:OU_injectivity} below, applied on $[s,T]$, shows that Equation (\ref{eq:zero_condition}) amounts to $g \overset{L^2([s,T],dt)}{=} 0$ and thus
\begin{equation}\label{eq:affine_dependancy}
\sum\limits_{i \in I} \alpha_i f^X_i \overset{L^2([s,T],dt)}{=} 0.
\end{equation}
\par \noindent Again, as $(f_i^X)_{i \in I}$ are linearly independent, we have $\alpha_1 = \cdots = \alpha_n = 0$. \emph{Hence the (\ref{eq:H_hypothesis}) hypothesis is fulfilled in the case of K-L generalized bridges of the Ornstein-Uhlenbeck processes.}
\vspace{2mm}
\begin{prop}[Injectivity of the Wiener integral related to centered Ornstein-Uhlenbeck processes]\label{prop:OU_injectivity}
\par Let $G$ be an Ornstein-Uhlenbeck process defined on $[0,T]$ by the SDE
$$
dG_t = -\theta G_t dt + \sigma dW_t \quad \textnormal{with } \sigma > 0 \textnormal{ and } \theta >0,
$$
where $W$ is a standard Brownian motion and $G_0 \simdist \mathcal{N}\left(0,\sigma_0^2\right)$ is independent of $W$. 
\par \noindent If $g \in L^2([0,T])$, then we have
$$
\int_0^T g(s) dG_s =0 \quad \Leftrightarrow \quad g \overset{L^2([0,T])}{=} 0. 
$$
\end{prop}
\par \noindent \textbf{Proof:} The solution of the Ornstein-Uhlenbeck SDE is $G_t = \!\!\!\!\! \underbrace{G_0 e^{- \theta t}}_{\textnormal{independent of } W}\!\!\! \overset{\independent }{+} \underbrace{\int_0^t \sigma e^{\theta (s-t)} dW_s}_{=: G_t^0}$. Hence, we have
$$
\int_0^T g(s) dG_s = -\theta G_0 \int_0^T g(s) e^{-\theta s}ds \overset{\independent }{+} \int_0^T g(s) dG_s^0.
$$
Thus, by independence, if $\int_0^T g(s) dG_s = 0$ then $\int_0^T g(s) dG_s^0 = 0$. This means that we only have to prove the proposition in the case of an Ornstein-Uhlenbeck process starting from $0$.
\vspace{2mm}
\par We now assume that $\sigma_0^2 = 0$ and we temporarily make the additional assumption that $\theta T < \frac{4}{3}$. If $g \in L^2([0,T])$ and $\int_0^T g(s) dG_s =0$, then $\theta \int_0^T g(s) G_s ds = \sigma \int_0^T g(s) dW_s$, and thus, if $\Gamma^{OU}$ denotes the covariance function of $G$, 
\begin{equation}\label{eq:variance_equality}
\theta^2\int_0^T\int_0^T g(s) g(t) \Gamma^{OU}(s,t) ds dt = \sigma^2 \int_0^T g(s)^2 ds.
\end{equation}
\par \noindent Applying Schwarz's inequality twice, we get
$$
\int_0^T\int_0^T g(s) g(t) \Gamma^{OU}(s,t) ds dt \leq \int_0^T g(s)^2 ds \sqrt{\int_0^T \int_0^T \left(\Gamma^{OU}(s,t)\right)^2 ds dt}.
$$
Hence, provided that 
\begin{equation}\label{eq:OU_covariance_inequality}
\int_0^T \int_0^T \left(\Gamma^{OU}(s,t)\right)^2 ds dt < \frac{\sigma^4}{\theta^4},
\end{equation}
Equality (\ref{eq:variance_equality}) implies $\int_0^T g(s)^2 ds =0$. 
\par \noindent Now, we come to the proof of Inequality (\ref{eq:OU_covariance_inequality}). The covariance function of the Ornstein-Uhlenbeck process starting from $0$ writes
$$
\Gamma^{OU}(s,t) = \frac{\sigma^2}{2\theta} e^{-\theta(s+t)} \left( e^{2 \theta \min(s,t)}-1\right).
$$
For $t \in [0,T]$, we have $\int_0^T \left(\Gamma^{OU}(s,t)\right)^2 ds = \frac{\sigma^4}{8\theta^3}\left(2 - 4e^{-2\theta t}\theta t - e^{-2 \theta (T-t)} - 2e^{-2 \theta t} +2 e^{-2 \theta T} - e^{-2 \theta (T+t)}\right)$, and thus
$$
\int_0^T \int_0^T \left( \Gamma^{OU}(s,t) \right)^2 ds dt = \frac{\sigma^2}{16\theta^4} \left( -5 +4\theta T + 8 \theta T e^{-2 \theta T} + 4e^{-2 \theta T} + e^{-4 \theta T }\right).
$$
\par \noindent Consequently, the function $\phi$ defined by $\phi(\theta) := \int_0^T \int_0^T \left(\Gamma^{OU}(s,t)\right)^2 ds dt - \frac{\sigma^4}{\theta^4}$ writes
$$
\phi(\theta) = \frac{1}{16} \frac{\sigma^4}{\theta^4} \left(-21+4 \theta T+8 \theta e^{-2 \theta T} T+4 e^{-2 \theta T}+e^{-4 \theta T}\right).
$$
\par \noindent Thus $\phi(\theta) < -16+12 \theta T$ which leads to Inequality (\ref{eq:OU_covariance_inequality}) thanks to the fact that $\theta T < \frac{4}{3}$.
\vspace{3mm}
\par We now come back to the general case where we might have $\theta T \geq \frac{4}{3}$. If this is the case, let us consider $\widetilde{T} := T-\frac{1}{\theta}$, so that $\theta \left(T-\widetilde{T}\right) < \frac{4}{3}$. For $t \in \left[\widetilde{T},T\right]$, we have 
$$
G_t = \underbrace{G_{\widetilde{T}} e^{- \theta \left(t-\widetilde{T}\right)}}_{\textnormal{independent of } (W_s)_{s\in [\widetilde{T},T]}} \overset{\independent }{+} \quad \underbrace{\int_{\widetilde{T}}^t \sigma e^{\theta (s-t)} dW_s}_{=: \widetilde{G}_t^0}.
$$
The so-defined process $\left(\widetilde{G}_t^0\right)_{t \in \left[\widetilde{T},T\right]}$ is a centered Ornstein-Uhenbeck process starting from $0$ and satisfying the same SDE as $G$. Hence, by independence, if $\int_0^T g(s) dG_s = 0$, then $\int_{\widetilde{T}}^T g(s) d\widetilde{G}_s^0 = 0$.
\par As $\theta \left(T-\widetilde{T}\right) < \frac{4}{3}$, we can apply the result to $\left(\widetilde{G}_t^0\right)_{t \in \left[\widetilde{T},T\right]}$ so that $g_{|\left[\widetilde{T},T\right]} \overset{L^2([\widetilde{T},T])}{=} 0$. If $\widetilde{T}\theta < \frac{4}{3}$, we then have $g \overset{L^2([0,T])}{=} 0$. If it is not the case, we use the same method by using the decomposition of $\left[0,\widetilde{T}\right]$ into $\left[0,\widetilde{T} - \frac{1}{\theta}\right]$ and $\left[\widetilde{T} - \frac{1}{\theta},\widetilde{T}\right]$ and so on. An easy induction finally shows that $g \overset{L^2([0,T])}{=} 0$.
\vspace{2mm}
\par \noindent The converse is obvious. \myqed
\vspace{2mm}
\par \textbf{The case of a more general Gaussian semimartingale}
\par \noindent In Appendix \ref{sec:wiener_appendix}, we investigate the problem for more general Gaussian semimartingales. As we have seen in the case of the Ornstein-Uhlenbeck process, if functions $\left(f_i\right)_{i \in I}$ are linearly independent in $L^2([s,T],d \langle X \rangle)$ for $s \in [0,T)$, the (\ref{eq:H_hypothesis}) hypothesis comes to the injectivity of the Wiener integral with respect to $X$ on $\sspan(f_i)_{i \in I}$ (on interval $[s,T]$).
\section{K-L generalized bridges and partial functional quantization}\label{sec:KLgeneralizedBridges}
\par \noindent We keep the notations and assumptions of Section \ref{sec:KL_generalized_bridge}. As we have seen, Equation (\ref{eq:conditioning_KL}) decomposes the process $X$ as the sum of a linear combination of the Karhunen-Loève coordinates $Y := (Y_i)_{i \in I}$ and an independent remainder term. We now consider $\widehat{Y}^\Gamma$ a stationary Voronoi $N$-quantization of $Y$. $\widehat{Y}^\Gamma$ can be written as a nearest neighbor projection of $Y$ on a finite codebook $\Gamma = (\gamma_1,\cdots,\gamma_N)$.
$$
\widehat{Y}^\Gamma = \Proj_\Gamma(Y), \hspace{5mm} \textnormal{where }\Proj_\Gamma \textnormal{ is a nearest neighbor projection on } \Gamma.
$$
For example, $\widehat{Y}^\Gamma$ can be a stationary product quantization or an optimal quadratic quantization of $Y$. We now define the stochastic process $\widetilde{X}^{I,\Gamma}$ by replacing $Y$ by $\widehat{Y}^\Gamma$ in the decomposition (\ref{eq:conditioning_KL}). We denote $\widetilde{X}^{I,\Gamma} = \Proj_{I,\Gamma}(X)$.
$$
\widetilde{X}^{I,\Gamma} = \sum\limits_{i \in I} \widehat{Y}^\Gamma_i e_i^X {\overset{\independent }{+}} \sum\limits_{i \in \N^* \backslash I} \sqrt{\lambda_i^X} \xi_i e_i^X.
$$
The conditional distribution of $\widetilde{X}^{I,\Gamma}$ given that $Y$ falls in the Voronoi cell of $\gamma_k$ is the probability distribution of the K-L generalized bridge with end-point $\gamma_k$. In other words, we have quantized the Karhunen-Loève coordinates of $X$ corresponding to $i \in I$, and not the other ones. 
\par \noindent The so-defined process $\widetilde{X}^{I,\Gamma}$ is called a \emph{partial functional quantization of $X$}. 
\subsection{Partial functional quantization of stochastic differential equations}\label{sec:SDE_partial_functional_quantization}
\par \noindent Let $X$ be a continuous centered Gaussian semimartingale on $[0,T]$ with $X_0 = 0$. We consider the SDE
\begin{equation}\label{eq:SDE_stratonovich_partial_quant}
dS_t = b(t,S_t) dt+ \sigma(t,S_t) dX_t, \hspace{5mm} S_0=x \in \R, \hspace{5mm} \textnormal{and } t\in [0,T],
\end{equation}
where $b(t,x)$ and $\sigma(t,x)$ are Borel functions, Lipschitz continuous with respect to $x$ uniformly in $t$, $\sigma$ and $b(\cdot,0)$ are bounded. This SDE admits a unique strong solution $S$. 
\vspace{2mm}
\par \noindent The conditional distribution given that $Y_i = y_i$ for $i \in I$ of $S$ is the strong solution of the stochastic differential equation $dS_t = b(t,S_t) dt+ \sigma(t,S_t) dX^{I,\overline{y}}_t$, with $S_0=x \in \R$, and for $t\in [0,T]$, where $X^{I,\overline{y}}_t$ is the corresponding K-L generalized bridge. 
\vspace{2mm}
\par \noindent Under the (\ref{eq:H_hypothesis}) hypothesis, this suggests to define the partial quantization of $S$ from a partial quantization $\widetilde{X}^{I,\Gamma}$ of $X$ by replacing $X$ by $\widetilde{X}^{I,\Gamma}$ in the SDE (\ref{eq:SDE_stratonovich_partial_quant}). We define the \emph{partial quantization} $\widetilde{S}^{I,\Gamma}$ as the process whose conditional distribution given that $Y$ falls in the Voronoi cell of $\gamma_k$ is the strong solution of the same SDE where $X$ is replaced by the K-L generalized bridge with end-point $\gamma_k$. We write
\begin{equation}\label{eq:diffusion_partial_quantization_definition}
d\widetilde{S}^{I,\Gamma}_t = b\left(t,\widetilde{S}^{I,\Gamma}_t\right)dt + \sigma\left(t,\widetilde{S}^{I,\Gamma}_t\right) d \widetilde{X}_t^{I,\Gamma}.
\end{equation}
\begin{remark}
The SDE is written in the Itô sens unlike in the previous works on full functional quantization \cite{PagesSellamiSDE,PagesPrintemsFunctional4} where the SDE was written in the Stratonovich sense. 
\end{remark}
\par \noindent Here, the set $I$ of quantized Karhunen-Loève coordinates does not depend on the quantization level, while in the case of full functional quantization, optimality is reached by adapting the quantization dimension. The optimal quantization dimension (or critical dimension) has been thoroughly investigated in \cite{LuschgyPagesFunctional3, LuschgyPagesFunctional2} and is shown to be asymptotically equivalent to the logarithm of the quantization level when in goes to infinity, in the cases of Brownian motion, Brownian bridge and Ornstein-Uhlenbeck processes. 
\subsection{Convergence of partially quantized SDEs}
\par We start by stating some useful inequalities for the sequel. Then we recall the so-called Zador's theorem which will be used in the proof of the $a.s.$ convergence of partially quantized SDEs. 
\begin{lemm}[Gronwall inequality for locally finite measures]\label{lemm:gronwall}
\par \noindent Consider $\mathcal{I}$ an interval of the form $[a,b)$ or $[a,b]$ with $a<b$ or $[a,\infty)$. Let $\mu$ be a locally finite measure on the Borel $\sigma$-algebra of $\mathcal{I}$. We consider $u$ a measurable function defined on $\mathcal{I}$ such that for all $t \in \mathcal{I}$, $\int_a^t |u(s)|\mu(ds) < +\infty$. We assume that there exists a Borel function $\psi$ on $\mathcal{I}$ such that
$$
u(t) \leq \psi(t) + \int_{[a,t)} u(s) \mu(ds), \quad \forall t \in \mathcal{I}.
$$
\begin{tabular}{l|l}
If 	& either $\psi$ is non-negative,\\
	& or $t \mapsto \mu([a,t))$ is continuous on $\mathcal{I}$ and for all $t \in \mathcal{I}$, $\int_a^t |\psi(s)|\mu(ds) < \infty$, \\
\end{tabular}
\vspace{2mm}
\par \noindent then $u$ satisfies the Gronwall inequality.
$$
u(t) \leq \psi(t) + \int_{[a,t)} \psi(s) \exp(\mu([s,t))) \mu(ds).
$$
\end{lemm}
\par \noindent A proof of this result is available in \cite[Appendix $5.1$]{EthierKurtzMarkovProcessses}.
\begin{lemm}[A Gronwall-like inequality in the non-decreasing case]\label{lemm:gronwall_like}
\par \noindent Consider $\mathcal{I}$ an interval of the form $[a,b)$ or $[a,b]$ with $a<b$ or $[a,\infty)$. Let $\mu$ be a locally finite measure on the Borel $\sigma$-algebra of $\mathcal{I}$. We consider $u$ a measurable \emph{non-decreasing} function defined on $\mathcal{I}$ such that for all $t \in \mathcal{I}$, $\int_a^t |u(s)|\mu(ds) < +\infty$. We assume that there exists a Borel function $\psi$ on $\mathcal{I}$, and two non-negative constants $(A,B) \in \R_+^2$ such that
\begin{equation}\label{eq:gronwall_like_hypothesis}
u(t) \leq \psi(t) + A \int_{[a,t)} u(s) \mu(ds) + B \sqrt{\int_{[a,t)} u(s)^2 \mu(ds)}, \quad \forall t \in \mathcal{I}.
\end{equation}
\begin{tabular}{l|l}
If \quad	& either $\psi$ is non-negative,\\
			& or $t \mapsto \mu([a,t))$ is continuous on $\mathcal{I}$ and for all $t \in \mathcal{I}$, $\int_a^t |\psi(s)|\mu(ds) < \infty$, \\
\end{tabular}
\vspace{2mm}
\par \noindent then $u$ satisfies the following Gronwall inequality.
$$
u(t) \leq 2\psi(t) + 2\left(2A+B^2\right) \int_{[a,t)} \psi(s) \exp\left(\left(2A+B^2\right) \mu([s,t))\right) \mu(ds).
$$
\end{lemm}
\par \noindent \textbf{Proof:} Using that for $(x,y)\in \R_+^2$, $\sqrt{xy} \leq \frac{1}{2}\left(\frac{x}{B}+B y \right)$, we have 
$$
\left(\int_{[a,t)} u(s)^2 \mu(ds) \right)^{\frac{1}{2}} \leq \left( u(t) \int_{[a,t)} u(s) \mu(ds) \right)^{\frac{1}{2}} \leq \frac{u(t)}{2B} + \frac{B}{2} \int_{[a,t)} u(s) \mu(ds).
$$
Plugging this in Inequality (\ref{eq:gronwall_like_hypothesis}) yields
$$
u(t) \leq 2\psi(t) + \left(2A+B^2\right) \int_{[a,t)} u(s) \mu(ds).
$$
Applying the regular Gronwall's inequality (Lemma \ref{lemm:gronwall}) yields the announced result. \myqed
\begin{theo}[Zador, Bucklew, Wise, Graf, Luschgy, Pagès]\label{thm:Zador}
\par \noindent Consider $r>0$ and $X$ be a $\R^d$-valued random variable such that $X \in L^{r+\eta}$ for some $\eta > 0$. We denote by $\mathcal{E}_{N,r}(X)$ the $L^r$ optimal quantization error of level $N$ for $X$, $\mathcal{E}_{N,r}(X) := \min\left\{ \left\|X - Y\right\|_r, |Y(\Omega)| \leq N \right\}$.
\begin{enumerate}
\item \emph{(Sharp rate)} Let $\PP_X(d\xi) = \phi(\xi) d \xi + \nu(d\xi)$ be the Radon-Nikodym decomposition of the probability distribution of $X$. ($\nu$ and the Lebesgue's measure are singular). Then if $\phi \not\equiv 0$,
$$
\mathcal{E}_{N,r}(X) \underset{N \to \infty}{\sim} \widetilde{J}_{r,d} \times \left( \int_{\R^d} \phi^{\frac{d}{d+r}}(u) du \right)^{\frac{1}{d}+\frac{1}{r}} \times N^{-\frac{1}{d}},
$$
where $\widetilde{J}_{r,d} \in (0, \infty)$.
\item \emph{(Non-asymptotic upper bound)} There exists $C_{d,r,\eta} \in (0,\infty)$ such that, for every $\R^d$-valued random vector $X$,
$$
\forall N \geq 1, \quad \mathcal{E}_{N,r}(X) \leq C_{d,r,\eta} \|X\|_{r + \eta} N^{-\frac{1}{d}}.
$$
\end{enumerate}
\end{theo}
\par \noindent The first statement of the theorem was first established for probability distributions with compact support by Zador \cite{ZadorAsymptoticError}, and extended by Bucklew and Wise to general probability distributions on $\R^d$ \cite{BuckleyWise}. The first mathematically rigorous proof can be found in \cite{GrafLushgyMonograf}. The proof of the second statement is available in \cite{FunctionalQuantizationLevy}. 
\par \noindent The real constant $\widetilde{J}_{r,d}$ corresponds to the case of the uniform probability distribution over the unit hypercube $[0,1]^d$. We have $\widetilde{J}_{r,1} = \frac{1}{2}(r+1)^{-\frac{1}{r}}$ and $\widetilde{J}_{2,2} = \sqrt{\frac{5}{18 \sqrt{2}}}$ (see \cite{GrafLushgyMonograf}).

\subsubsection{$L^p$ convergence of partially quantized SDEs}
\begin{lemm}[Generalized Minkowski inequality for locally finite measures]\label{lemm:generalized_minkowski}
\par \noindent Consider $\mathcal{I}$ an interval of the form $[a,b)$ or $[a,b]$ with $a < b$ or $[a,\infty)$. Let $\mu$ be a locally finite measure on the Borel $\sigma$-algebra of $\mathcal{I}$. Then for any non-negative bi-measurable process $X = (X_t)_{t \in \mathcal{I}}$ and every $p \in [1,\infty)$,
$$
\left\| \int_{\mathcal{I}} X_t \mu(dt) \right\|_p \leq \int_{\mathcal{I}} \|X_t\|_p \mu(dt).
$$
\end{lemm}
\begin{prop}[Burkholder-Davis-Gundy inequality]
\par \noindent For every $p\in (0,\infty)$, there exist two positive real constants $c_p^{BDG}$ and $C_p^{BDG}$ such that for every continuous local martingale $(X_t)_{t \in [0,T]}$ null at $0$,
$$
c_p^{BDG} \left\| \sqrt{\langle X \rangle_T} \right\|_p \leq \left\|\sup\limits_{s \in [0,T]} |X_s| \right\|_p \leq C_p^{BDG} \left\| \sqrt{\langle X \rangle_T} \right\|_p.
$$
\end{prop}
\par \noindent We refer to \cite{RevuzYor} for a detailed proof. 
\begin{prop}[$L^p$ inequality]
\par Let $G$ be a standard Gaussian random variable valued in $\R$. There exists a constant $C_p > 0$ such that for every $M>1$
$$
\sqrt{\frac{2}{\pi}} M^{p-1} \exp \left( -\frac{M^2}{2} \right) \leq \E\left[ |G|^p \1_{|G|>M} \right] \leq C_p M^{p-1} \exp\left(-\frac{M^2}{2} \right).
$$
Consequently
$$
\left(\sqrt{\frac{2}{\pi}}\right)^{\frac{1}{p}} M^{\frac{1}{q}} \exp\left(-\frac{M^2}{2p} \right) \leq \left\|G \1_{|G|>M} \right\|_p \leq (C_p)^{1/p} M^{\frac{1}{q}} \exp\left(-\frac{M^2}{2p} \right),
$$
where $q$ is the conjugate exponent of $p$. 
\end{prop}
\begin{prop}[The non-standard case and $L^p$ reverse inequality]\label{prop:Lp_non_standard_inverse}
\par \noindent If $H := \sigma G$ has a variance of $\sigma^2$, we obtain
\begin{equation}\label{eq:Lp_general_inequality}
\begin{array}{lll}
\left\|H \1_{|H|>M} \right\|_p	&\leq	& \sigma \left\|G \1_{|G|>\frac{M}{\sigma}} \right\|_p = \sigma (C_p)^{1/p} \left(\frac{M}{\sigma}\right)^{\frac{1}{q}} \exp\left(-\frac{M^2}{2p\sigma^2} \right),\\
								&		& = \underbrace{\sigma^{\frac{1}{p}} (C_p)^{1/p} M^{\frac{1}{q}} \exp\left(-\frac{M^2}{2p\sigma^2} \right)}_{=: \eta_M}.
\end{array}
\end{equation}
\par \noindent Conversely, for some fixed $\eta>0$, and if $M>1$, we have
\begin{equation}\label{eq:Lp_reverse_inequality}
M \geq \underbrace{\sqrt{-\sigma^2 (p-1) \mathcal{W}_{-1}\left(-\frac{q \eta^{2q}}{p \sigma^2 \left( C_p^{2q/p} \sigma^{2q/p}\right)}\right)}}_{=: M_\eta} \Rightarrow \eta_M \leq \eta \, 
\end{equation}
where $\mathcal{W}_{-1}$ is the secondary branch of the Lambert $\mathcal{W}$ function. For more details on the Lambert $\mathcal{W}$ function, we refer to \cite{CorlessKnuthLambert}.
\end{prop}
\begin{theo}[$L^p$ quantization of partially quantized SDEs]\label{thm:p_partial_quantization_convergence}
\par \noindent Let $X$ be a continuous centered Gaussian martingale on $[0,T]$ with $X_0 = 0$. Let $S$ be the strong solution of the SDE
$$
dS_t = b(t,S_t) dt + \sigma(t,S_t) dX_t, \quad S_0=x,
$$
where $b(t,x)$ and $\sigma(t,x)$ are Borel functions, Lipschitz continuous with respect to $x$ uniformly in $t$, $\sigma$ and $b(\cdot,0)$ are bounded.
\par \noindent We consider $\widetilde{X}^{I,\Gamma}$ a stationary partial functional quantization of $X$ and $\widetilde{S}^{I,\Gamma}$ the corresponding partial functional quantization of $S$, \textit{i.e.} the strong solutions of
$$
d\widetilde{S}^{I,\Gamma}_t = b\left(t,\widetilde{S}^{I,\Gamma}_t\right) dt + \sigma\left(t,\widetilde{S}^{I,\Gamma}_t\right) d\widetilde{X}^{I,\Gamma}_t, \quad \widetilde{S}^{I,\Gamma}_0=x.
$$
Then, for every $p \in (0,\infty)$, $\varepsilon >0$ and $t \in [0,T)$, there exist three positive constants $C_{X,\varepsilon,I}$, $A_{X,\varepsilon,I}$ and $B_{X,\varepsilon,I}$ such that
\begin{equation}\label{eq:diffusion_partial_quantization_p_err}
\left\| \sup\limits_{v \in [0,t]} \left|S_v - \widetilde{S}^{I,\Gamma}_v \right| \right\|_p \leq C_{X,\varepsilon,I} \exp\left(A_{X,\varepsilon,I} \sqrt{-\mathcal{W}_{-1}\left(-\frac{ \left\|Y-\widehat{Y}^\Gamma \right\|_{p+\varepsilon}^{2q}}{B_{X,\varepsilon,I}} \right)} \right) \left\|Y-\widehat{Y}^\Gamma \right\|_{p+\varepsilon},
\end{equation}
where $q$ is the conjugate exponent of $p$, where $Y$ is defined from $X$ by Equation (\ref{eq:conditioning_KL}) and $\widehat{Y}^{\Gamma}$ is the nearest neighbor projection on $\Gamma$.
\end{theo}
\begin{remark}
Using that $\mathcal{W}_{-1}(-x) \underset{x \to 0_+}{\sim} \ln(x)$, we can see that the right-hand term in Equation (\ref{eq:diffusion_partial_quantization_p_err}) goes to $0$ as the quantization error $\left\|Y-\widehat{Y}^\Gamma \right\|_{p+\varepsilon}$ goes to $0$.
\end{remark}
\par \noindent \textbf{Proof:} We decompose the process $X$ into $X_t = \sum\limits_{i\in I} Y_i e_i^X(t) + X_t^{I,\overline{0}}$ and $\widetilde{X}^{I,\Gamma}$ into $\widetilde{X}^{I,\Gamma}_t = \sum\limits_{i \in I} \widehat{Y}_i^\Gamma e_i^X(t) + X_t^{I,\overline{0}}$, where $\widehat{Y}^{\Gamma}$ is the nearest neighbor projection of $Y$ on $\Gamma$. For some $k \in \{1, \cdots,N\}$, conditionally to $\widehat{Y}^\Gamma = \gamma_k$, we have
\begin{multline*}
S_t-\widetilde{S}^{I,\Gamma}_t = \int_0^t \left(b(u,S_u) - b\left(u,\widetilde{S}^{I,\Gamma}_t\right)\right)du + \sum\limits_{i \in I} \int_0^t \left(\sigma(u,S_u) - \sigma\left(u,\widetilde{S}^{I,\Gamma}_u\right) \right) \widehat{Y}^\Gamma_i d e_i^X(u) \\
+\sum\limits_{i \in I} \int_0^t \left(Y_i - \widehat{Y}^\Gamma_i\right) \sigma(u,S_u) de_i^X(u) + \int_0^t \left(\sigma(u,S_u) - \sigma\left(u,\widetilde{S}^{I,\Gamma}_u\right) \right) G_u d\langle X \rangle_u\\
+ \int_0^t \left( \sigma(u,S_u) - \sigma\left(u,\widetilde{S}^{I,\Gamma}_u\right) \right) d \widetilde{M}_u.
\end{multline*}
This gives (conditionally to $\widehat{Y}^\Gamma = \gamma_k$)
\begin{multline*}
\left|S_t - \widetilde{S}^{I,\Gamma}_t\right| \leq [b]_{\textnormal{Lip}} \int_0^t \left|S_u - \widetilde{S}^{I,\Gamma}_u \right|du + [\sigma]_{\textnormal{Lip}} |I| \max\limits_{\substack{i \in I\\ u \in [0,T]}} \left| \left(e_i^X\right)'(u)\right| \left(\max\limits_{i \in I} \left|\widehat{Y}^\Gamma_i \right| \right) \int_0^t \left| S_u - \widetilde{S}^{I,\Gamma}_u\right| du\\
+ [\sigma]_{\max} |I| \max\limits_{\substack{i \in I \\ u \in [0,T]}} \left| \left(e_i^X\right)'(u)\right| T \sum\limits_{i \in I} \left|Y_i - \widehat{Y}^\Gamma_i \right| + \left|\int_0^t \left( \sigma(u,S_u) - \sigma \left(u,\widetilde{S}^{I,\Gamma}_u\right)\right) G_u d \langle X \rangle_u \right|\\
 +\left| \int_0^t \left( \sigma(u,S_u) - \sigma \left(u,\widetilde{S}^{I,\Gamma}_u)\right)\right) d \widetilde{M}_u \right|.
\end{multline*}
As a consequence, conditionally to $\widehat{Y}^\Gamma = \gamma_k$,
\begin{multline*}
\max\limits_{v \in [0,t]} \left|S_v -\widetilde{S}^{I,\Gamma}_v\right| \leq [b]_{\textnormal{Lip}} \int_0^t \max_{v \in [0,u]} \left|S_v - \widetilde{S}^{I,\Gamma}_v\right| du\\
+ [\sigma]_{\textnormal{Lip}} |I| \max\limits_{\substack{i \in I\\ u \in [0,T]}} \left| \left( e_i^X \right)'(u)\right| \left(\max\limits_{i \in I} \left|\widehat{Y}^\Gamma_i \right| \right) \int_0^t \max\limits_{v \in [0,u]} \left| S_v -\widetilde{S}^{I,\Gamma}_v \right| du \\
+ [\sigma]_{\max} |I| \max\limits_{\substack{i \in I \\ u \in [0,T]}} \left|\left(e_i^X\right)'(u)\right| T \sum\limits_{i \in I} \left|Y_i - \widehat{Y}^\Gamma_i\right| + \max\limits_{v \in [0,t]} \left|\int_0^v \left(\sigma(u,S_u) - \sigma\left(u,\widetilde{S}^{I,\Gamma}_u \right)\right) G_u d \langle X \rangle_u \right|\\
+ \max\limits_{v \in [0,t]} \left|\int_0^v \left(\sigma(u,S_u) - \sigma\left(u,\widetilde{S}^{I,\Gamma}_u\right)\right) d \widetilde{M}_u \right|.
\end{multline*}
\par \noindent To shorten the notations, we denote, for a random variable $V$ and a non-negligible event $A$, $\|V\|_{p,A} := \E\left[V^p \middle| A \right]^{1/p}$. Hence, using the Minkowski inequality and the generalized Minkowski inequality for locally finite measures (Lemma \ref{lemm:generalized_minkowski}), we get
\begin{multline*}
\left\| \max\limits_{v \in [0,t]} \left|S_v -\widetilde{S}^{I,\Gamma}_v\right| \right\|_{p,\{\widehat{Y}^\Gamma = \gamma_k \}} \leq [b]_{\textnormal{Lip}} \int_0^t \left\|\max_{v \in [0,u]} \left|S_v - \widetilde{S}^{I,\Gamma}_v\right| \right\|_{p,\{\widehat{Y}^\Gamma = \gamma_k \}} du\\
+ [\sigma]_{\textnormal{Lip}} |I| \max\limits_{\substack{i \in I \\ u \in [0,T]}} \left|\left(e_i^X\right)'(u) \right| \left( \max\limits_{i\in I} \left| \widehat{Y}_i^\Gamma \right| \right) \int_0^t \left\|\max_{v \in [0,u]} \left|S_v - \widetilde{S}^{I,\Gamma}_v\right| \right\|_{p,\{\widehat{Y}^\Gamma = \gamma_k \}} du\\
+ [\sigma]_{\textnormal{Lip}} |I| \max\limits_{\substack{i \in I \\ u \in [0,T]}} \left|\left(e_i^X\right)'(u) \right| T \left\| \sum\limits_{i \in I} \left| Y_i - \widehat{Y}^\Gamma_i \right| \right\|_{p,\{\widehat{Y}^\Gamma = \gamma_k \}}
+ \left\| \max\limits_{v \in [0,t]} \left|\int_0^v \left(\sigma(u,S_u) - \sigma\left(u,\widetilde{S}^{I,\Gamma}_u \right)\right) G_u d \langle X \rangle_u \right| \right\|_{p,\{\widehat{Y}^\Gamma = \gamma_k \}}\\
+ \left\|\max\limits_{v \in [0,t]} \left|\int_0^v \left(\sigma(u,S_u) - \sigma\left(u,\widetilde{S}^{I,\Gamma}_u\right)\right) d \widetilde{M}_u \right| \right\|_{p,\{\widehat{Y}^\Gamma = \gamma_k \}}.
\end{multline*}
\par \noindent Now, from the Burkholder-Davis-Gundy inequality,
\begin{multline}\label{eq:Lp_intermediate_inequality}
\left\| \max\limits_{v \in [0,t]} \left|S_v -\widetilde{S}^{I,\Gamma}_v\right| \right\|_{p,\{\widehat{Y}^\Gamma = \gamma_k \}} \leq [b]_{\textnormal{Lip}} \int_0^t \left\|\max_{v \in [0,u]} \left|S_v - \widetilde{S}^{I,\Gamma}_v\right| \right\|_{p,\{\widehat{Y}^\Gamma = \gamma_k \}} du\\
+ [\sigma]_{\textnormal{Lip}} |I| \max\limits_{\substack{i \in I \\ u \in [0,T]}} \left|\left(e_i^X\right)'(u) \right| \left( \max\limits_{i\in I} \left| \widehat{Y}_i^\Gamma \right| \right) \int_0^t \left\|\max_{v \in [0,u]} \left|S_v - \widetilde{S}^{I,\Gamma}_v\right| \right\|_{p,\{\widehat{Y}^\Gamma = \gamma_k \}} du\\
+ [\sigma]_{\textnormal{Lip}} |I| \max\limits_{\substack{i \in I \\ u \in [0,T]}} \left|\left(e_i^X\right)'(u) \right| T \left\| \sum\limits_{i \in I} \left| Y_i - \widehat{Y}^\Gamma_i \right| \right\|_{p,\{\widehat{Y}^\Gamma = \gamma_k \}}\\
+ \left\|\int_0^t \left|\sigma(u,S_u) - \sigma\left(u,\widetilde{S}^{I,\Gamma}_u \right)\right| |G_u| d \langle X \rangle_u \right\|_{p,\{\widehat{Y}^\Gamma = \gamma_k \}}\\
+ C_p^{BDG} \left\| \sqrt{\int_0^t \left( \sigma(u,S_u) - \sigma\left(u,\widetilde{S}^{I,\Gamma}_u \right) \right)^2 d \langle X \rangle_u} \right\|_{p,\{\widehat{Y}^\Gamma = \gamma_k \}}.
\end{multline}
\par \noindent Now, from Schwarz's inequality
$$
\left\| \sum\limits_{i \in I} \left| Y_i - \widehat{Y}^\Gamma_i \right| \right\|_{p,\{\widehat{Y}^\Gamma = \gamma_k \}} \leq \left\| \sqrt{|I|} \sqrt{\sum\limits_{i \in I} \left| Y_i - \widehat{Y}^\Gamma_i \right|^2 } \right\|_{p,\{\widehat{Y}^\Gamma = \gamma_k \}} = \sqrt{|I|} \left\| Y - \widehat{Y}^\Gamma \right\|_{p,\{\widehat{Y}^\Gamma = \gamma_k \}}.
$$
\par \noindent From the generalized Minkowski inequality
$$
\begin{array}{lll}
\left\|\int_0^t \left|\sigma(u,S_u) - \sigma\left(u,\widetilde{S}^{I,\Gamma}_u \right)\right| |G_u| d \langle X \rangle_u \right\|_{p,\{\widehat{Y}^\Gamma = \gamma_k \}} \leq \int_0^t \left\| \left(\sigma(u,S_u) - \sigma\left(u,\widetilde{S}^{I,\Gamma}_u \right) \right) G_u \right\|_{p,\{\widehat{Y}^\Gamma = \gamma_k \}} d \langle X \rangle_u\\
\hfill = \int_0^t \left\|\left(\sigma(u,S_u) - \sigma\left(u,\widetilde{S}^{I,\Gamma}_u \right) \right) G_u \1_{|G_u| \geq M} + \left(\sigma(u,S_u) - \sigma\left(u,\widetilde{S}^{I,\Gamma}_u \right) \right) G_u \1_{|G_u| \leq M} \right\|_{p,\{\widehat{Y}^\Gamma = \gamma_k \}} d\langle X \rangle_u\\
\qquad \leq \int_0^t \left\|\left(\sigma(u,S_u) - \sigma\left(u,\widetilde{S}^{I,\Gamma}_u \right) \right) G_u \1_{|G_u| \geq M}\right\|_{p,\{\widehat{Y}^\Gamma = \gamma_k \}} d\langle X \rangle_u \\
\hfill + \int_0^t \left\|\left(\sigma(u,S_u) - \sigma\left(u,\widetilde{S}^{I,\Gamma}_u \right) \right) G_u \1_{|G_u| \leq M} \right\|_{p,\{\widehat{Y}^\Gamma = \gamma_k \}} d\langle X \rangle_u\\
\qquad \leq 2 [\sigma]_{\max} \int_0^t \left\| G_u \1_{|G_u|\geq M} \right\|_{p,\{\widehat{Y}^\Gamma = \gamma_k \}} d\langle X \rangle _u + M [\sigma]_{\textnormal{Lip}} \int_0^t \left\| S_u - \widetilde{S}_u^{I,\Gamma}\right\|_{p,\{\widehat{Y}^\Gamma = \gamma_k \}} d\langle X \rangle_u.
\end{array}
$$
\par \noindent We obtain, thanks to Proposition \ref{prop:Lp_non_standard_inverse}
$$
\begin{array}{llll}
\left\|\int_0^t \left|\sigma(u,S_u) - \sigma\left(u,\widetilde{S}^{I,\Gamma}_u \right)\right| |G_u| d \langle X \rangle_u \right\|_{p,\{\widehat{Y}^\Gamma = \gamma_k \}}\\
\leq \underbrace{2 [\sigma]_{\max} \langle X \rangle_t (C_p)^{1/p} v_t^{\frac{1}{p}} M^{\frac{1}{q}} \exp\left(-\frac{M^2}{2p v_t^2}\right)}_{=: \eta_M}+ M [\sigma]_{\textnormal{Lip}} \int_0^t \left\| S_u - \widetilde{S}_u^{I,\Gamma}\right\|_{p,\{\widehat{Y}^\Gamma = \gamma_k \}} d\langle X \rangle_u,
\end{array}
$$
where $v_t^2 = \max\limits_{u \in [0,t]}\left(\var(G_u)\right)$. Moreover 
$$
\left\| \sqrt{\int_0^t \left( \sigma(u,S_u) - \sigma\left(u,\widetilde{S}^{I,\Gamma}_u \right) \right)^2 d \langle X \rangle_u} \right\|_{p,\{\widehat{Y}^\Gamma = \gamma_k \}} \leq \sqrt{\int_0^t \left\| \max\limits_{\substack{i \in I \\ v \in [0,u]}} \left| S_v - \widetilde{S}^{I,\Gamma}_v \right| \right\|_{p,\{\widehat{Y}^\Gamma = \gamma_k \}}^2 d \langle X \rangle_u}.
$$
\par \noindent Hence, Equation (\ref{eq:Lp_intermediate_inequality}) becomes
\begin{multline}
\left\| \max\limits_{v \in [0,t]} \left|S_v -\widetilde{S}^{I,\Gamma}_v\right| \right\|_{p,\{\widehat{Y}^\Gamma = \gamma_k \}} \leq \underbrace{[\sigma]_{\textnormal{Lip}} |I| \max\limits_{\substack{i \in I \\ u \in [0,T]}} \left|\left(e_i^X\right)'(u) \right| \sqrt{|I|}}_{=: A_i^X} \left\| Y - \widehat{Y}^\Gamma \right\|_{p,\{\widehat{Y}^\Gamma = \gamma_k \}} + \eta_M\\
+ [b]_{\textnormal{Lip}} \int_0^t \left\|\max_{v \in [0,u]} \left|S_v - \widetilde{S}^{I,\Gamma}_v\right| \right\|_{p,\{\widehat{Y}^\Gamma = \gamma_k \}} du\\
+ [\sigma]_{\textnormal{Lip}} |I| \max\limits_{\substack{i \in I \\ u \in [0,T]}} \left|\left(e_i^X\right)'(u) \right| \left( \max\limits_{i\in I} \left| \widehat{Y}_i^\Gamma \right| \right) \int_0^t \left\|\max_{v \in [0,u]} \left|S_v - \widetilde{S}^{I,\Gamma}_v\right| \right\|_{p,\{\widehat{Y}^\Gamma = \gamma_k \}} du\\
+ C_p^{BDG} \left( \int_0^t 2 \left\| \max\limits_{\substack{i \in I \\ v \in [0,u]}} \left| S_v - \widetilde{S}^{I,\Gamma}_v \right| \right\|_{p,\{\widehat{Y}^\Gamma = \gamma_k \}}^2 d \langle X \rangle_u \right)^{1/2}\\
+ \underbrace{M [\sigma]_{\textnormal{Lip}}}_{=: C^{X,M}} \int_0^t \left\| \max\limits_{v \in [0,u]} \left| S_v - \widetilde{S}_v^{I,\Gamma}\right| \right\|_{p,\{\widehat{Y}^\Gamma = \gamma_k \}} d\langle X \rangle_u.
\end{multline}
\par \noindent We can then apply the ``Gronwall-like'' lemma \ref{lemm:gronwall_like} for locally finite measures to the non-decreasing function 
$$
\left\| \sup\limits_{v\in [0,t]} \left|S_v-\widetilde{S}^{I,\Gamma}_v \right| \right\|_{p,\{\widehat{Y}^\Gamma = \gamma_k \}} = \E \left[ \sup\limits_{v\in [0,t]} \left|S_v-\widetilde{S}^{I,\Gamma}_v \right|^p \middle| \widehat{Y}^\Gamma = \gamma_k \right]^{1/p}
$$
and with the locally finite measure $\mu$ defined by $\mu(du) = du + d\langle X \rangle_u$, and we obtain
{\small
$$
\begin{array}{l}
\left\| \sup\limits_{v\in [0,t]} \left|S_v-\widetilde{S}^{I,\Gamma}_v \right| \right\|_{p,\{\widehat{Y}^\Gamma = \gamma_k \}} \leq \left( A_I^X \E\left[ \left|Y - \widehat{Y}^\Gamma \right|^p \middle| \widehat{Y}^\Gamma = \gamma_k \right]^{1/p} + \eta_M \right) \exp\left(\left( E_I^{X,\gamma_k} + C^{X,M} \right)\mu([0,t)) \right)\\
\hfill \leq \left( A_I^X \E\left[ \left|Y - \widehat{Y}^\Gamma \right|^p \middle| \widehat{Y}^\Gamma = \gamma_k \right]^{1/p} + \eta_M \right) \underbrace{\exp\left( E_I^{X,\gamma_k} \mu([0,t)) \right)}_{=: \phi(\gamma_k)} \exp\left(C^{X,M}\mu([0,t)) \right),
\end{array}
$$
}
where $E_I^{X,\gamma_k}$ is an affine function of $\max\limits_{i \in I} \left| \left(\gamma_k\right)_i \right|$. This yields
$$
\left\| \sup\limits_{v\in [0,t]} \left|S_v-\widetilde{S}^{I,\Gamma}_v \right| \right\|_p \leq \left( A_I^X \left\| \E\left[ \left|Y - \widehat{Y}^\Gamma \right|^p \middle| \widehat{Y}^\Gamma \right]^{1/p} \phi\left( \widehat{Y}^\Gamma \right) \right\|_p + \eta_M \left\| \phi\left( \widehat{Y}^\Gamma \right) \right\|_p \right) \exp\left(C^{X,M}\mu([0,t)) \right).
$$
\par \noindent Now, for $\varepsilon > 0$ and $\tilde{p}= 1 + \frac{\varepsilon}{p}$ and $\tilde{q} = \frac{\tilde{p}}{\tilde{p}-1} = 1 + \frac{p}{\varepsilon}$ the conjugate exponent of $\tilde{p}$, we have, thanks to Hölder's inequality
$$
\begin{array}{lll}
\E\left[ \phi\left( \widehat{Y}^\Gamma \right)^p \E\left[ \left|Y - \widehat{Y}^\Gamma \right|^p \middle| \widehat{Y}^\Gamma\right] \right] &\leq \left\| \phi\left(\widehat{Y}^\Gamma\right)^p \right\|_{\tilde{q}} \left\| \E\left[ \left|Y - \widehat{Y}^\Gamma \right|^p \middle|\widehat{Y}^\Gamma \right]\right\|_{\tilde{p}}\\
& \leq \left\| \phi\left(\widehat{Y}^\Gamma\right)^p \right\|_{\tilde{q}} \E\left[ \left|Y - \widehat{Y}^\Gamma \right|^{p+\varepsilon}\right]^{\frac{p}{p+\varepsilon}}.
\end{array}
$$
Hence,
$$
\left\| \E\left[ \left|Y - \widehat{Y}^\Gamma \right|^p \middle| \widehat{Y}^\Gamma \right]^{1/p} \phi\left( \widehat{Y}^\Gamma \right) \right\|_p \leq \left\| \phi\left(\widehat{Y}^\Gamma\right)^p \right\|_{\tilde{q}}^{1/p} \E\left[ \left|Y - \widehat{Y}^\Gamma \right|^{p+\varepsilon}\right]^{\frac{1}{p+\varepsilon}}.
$$
\par \noindent Now, as the so-defined function $\phi$ is convex and as $\widehat{Y}^\Gamma$ is a stationary quantizer of $Y$, we have thanks to Equation (\ref{eq:error_convex}), $\left\|\phi\left(\widehat{Y}^\Gamma\right)^p \right\|_{\tilde{q}} \leq \|\phi\left(Y\right)^p\|_{\tilde{q}}$ and $\left\|\phi\left(\widehat{Y}^\Gamma\right)\right\|_p \leq \|\phi\left(Y\right)\|_p$. 
\par \noindent If one sets $M= \sqrt{-v_t (p-1)\mathcal{W}_{-1}\left(-\frac{q\left\|Y-\widehat{Y}^\Gamma \right\|_{p+\varepsilon}^{2q}}{p v_t^2 C_p^{2q/p} v_t^{2q/p}} \right)}$, where $q$ is the conjugate exponent of $p$ and $\mathcal{W}_{-1}$ is the secondary branch of the Lambert $\mathcal{W}$ function, Proposition \ref{prop:Lp_non_standard_inverse} ensures that $\eta_M \leq \eta := \left\|Y-\widehat{Y}^\Gamma \right\|_{p + \varepsilon}$. We finally have the following error bound
{\small
$$
\left\| \sup\limits_{v \in [0,t]} \left|S_v - \widetilde{S}^{I,\Gamma}_v \right| \right\|_p \leq C_{X,\varepsilon,I} \exp\left([\sigma]_{\textnormal{Lip}} \sqrt{-v_t (p-1)\mathcal{W}_{-1}\left(-\frac{q\left\|Y-\widehat{Y}^\Gamma \right\|_{p+\varepsilon}^{2q}}{p v_t^2 C_p^{2q/p} v_t^{2q/p}} \right)} \right) \left\|Y-\widehat{Y}^\Gamma \right\|_{p+\varepsilon},
$$
}
which is the desired inequality \myqed
\begin{remark}[Without the stationarity property]\label{rq:without_stationarity}
\par \noindent The last step of the demonstration of Theorem \ref{thm:p_partial_quantization_convergence} (the use of Jensen's inequality) relies on the stationarity of the quantizer $\widehat{Y}$. Now, without this stationarity hypothesis and under the additional assumption that
\begin{equation}\tag{$\mathcal{A}$}\label{eq:not_empty}
\Gamma \cap B(0,1) \neq \emptyset,
\end{equation}
we have for every $i \in I$
$$
\left|\widehat{Y}_i \right| \leq \left|Y_i-\widehat{Y}_i \right| + |Y_i| \leq |Y_i| + \left|Y_i-\gamma^{k_0}_i\right| \leq 2|Y_i| + \left|\gamma^{k_0}_i\right| \leq 2|Y_i| + 1, \quad \textnormal{where } \gamma^{k_0} \in \Gamma \cap B(0,1).
$$
Hence
$$
\max\limits_{i\in I} \left|\widehat{Y}_i \right| \leq 2\max\limits_{i \in I} |Y_i| + 1.
$$
We notice that the function $\phi(x)$ defined in the demonstration of Theorem \ref{thm:p_partial_quantization_convergence} writes $\phi(x) = \psi(\max\limits_{i\in I} x_i)$ for some non-decreasing function $\psi$. This implies 
$$
\phi\left(\widehat{Y}\right) = \psi\left(\max\limits_{i\in I}\widehat{Y}_i\right) \leq \psi\left(\max\limits_{i\in I} \left(2|Y_i|+1 \right)\right) = \phi(2|Y|+1).
$$
Hence, we can obtain the same conclusion as in Theorem \ref{thm:p_partial_quantization_convergence}.
\end{remark}
\begin{coro}[$L^p$ convergence]\label{coro:lp_convergence}
\par \noindent With the same notations and hypothesis as in Theorem \ref{thm:p_partial_quantization_convergence}, consider $\left(\widetilde{X}^{I,\Gamma_n}\right)_{n \in \N}$ a sequence of partial functional quantizers of $X$ and $\left(\widetilde{S}^{I,\Gamma_n}\right)_{n \in \N}$ the corresponding sequence of partial quantizers of $S$. (For $n\in \N$, $\Gamma_n$ is assumed to have cardinal $n$.) 
\par \noindent If we make the additional assumption that the associated sequence of quantizers $\left(\widehat{Y}^{\Gamma_n}\right)_{n\in \N}$ is rate-optimal for the $L^{p+\varepsilon}$ convergence for some $\varepsilon>0$, then for every $t\in [0,T)$ we have
$$
\E\left[ \sup\limits_{u \in [0,t]} \left| S_u - \widetilde{S}_u^{I,\Gamma_n} \right|^p\right]^{1/p} = O \left(n^{-\frac{1}{|I|}} \right).
$$
\end{coro}
\par \noindent \textbf{Proof:} As $\left\|Y - \widehat{Y}^{\Gamma_n} \right\|_{p+\varepsilon} \underset{n\to \infty}{\to} 0$, we have $a.s.$ $d\left(\widehat{Y}^{\Gamma_n},Y \right) \underset{n\to \infty}{\to} 0$. Hence, there exists $N_0 \in \N$ such that for every $n\geq N_0$, $\Gamma_n$ satisfies Hypothesis (\ref{eq:not_empty}). From this observation, the result is straightforward consequence of the previous remark and Zador's theorem \ref{thm:Zador}, which defines the optimal convergence rate of a sequence of quantizers. \myqed
\subsubsection{The $a.s.$ convergence of partially quantized SDEs}
\begin{theo}[Almost sure convergence of partially quantized SDEs]\label{thm:as_partial_quantization_convergence}
\par Let $X$ be a continuous centered Gaussian martingale on $[0,T]$ with $X_0 = 0$. Let $S$ be the strong solution of the SDE
$$
dS_t = b(t,S_t) dt + \sigma(t,S_t) dX_t, \quad S_0=x,
$$
where $b(t,x)$ and $\sigma(t,x)$ are Borel functions, Lipschitz continuous with respect to $x$ uniformly in $t$, $\sigma$ and $b(\cdot,0)$ are bounded. 
\par \noindent We consider $\left(\widetilde{X}^{I,\Gamma_k}\right)_{k \in \N}$ a sequence of partial functional quantizers of $X$ and $\widetilde{S}^{I,\Gamma_n}$ the corresponding partial functional quantization of $S$, \textit{i.e.} the strong solutions of
$$
d\widetilde{S}^{I,\Gamma_n}_t = b\left(t,\widetilde{S}^{I,\Gamma_n}_t\right) dt + \sigma\left(t,\widetilde{S}^{I,\Gamma_n}_t\right) d\widetilde{X}^{I,\Gamma_n}_t, \quad \widetilde{S}^{I,\Gamma_n}_0=x.
$$
\par \noindent (For $n\in \N$, $\Gamma_n$ is assumed to have cardinal $n$.) We also assume that the sequence of partial quantizers of $X$ is rate-optimal for some $p > |I|$, \textit{i.e.} that there exists a constant $C$ such that
$$
\E\left[ \left|Y-\widehat{Y}^{\Gamma_n} \right|^p \right]^{1/p} \leq C n^{-\frac{1}{|I|}}
$$
for every $n \in \N^*$, where $Y$ is defined from $X$ by Equation (\ref{eq:conditioning_KL}) and $\widehat{Y}^{\Gamma}$ is the nearest neighbor projection on $\Gamma$. Then for every $t \in [0,T)$, $\widetilde{S}^{I,\Gamma_n}_t$ converges almost surely to $S_t$.
\end{theo}
\par \noindent \textbf{Proof:} From Corollary \ref{coro:lp_convergence}, if $t\in [0,T)$, there exist $r \in (|I|,p)$ and $N_0 \in \N$ such that for $n \geq N_0$,
$$
\E\left[ \sup\limits_{u \in [0,t]} \left| S_u - \widetilde{S}_u^{I,\Gamma_n} \right|^{r}\right]^{1/r} = O \left( n^{-\frac{1}{|I|}} \right).
$$
\par \noindent Hence, as $\frac{r}{|I|}>1$, Beppo-Levi's theorem for series with non-negative terms implies
$$
\E \left[ \sum\limits_{n \geq 1} \sup\limits_{u \in [0,t]} \left| S_u - \widetilde{S}_u^{I,\Gamma_n} \right|^{r} \right] < +\infty.
$$
Thus $\sum\limits_{n \geq 1} \sup\limits_{u \in [0,t]} \left| S_u - \widetilde{S}_u^{I,\Gamma_n} \right|^{r} < + \infty \ \PP - a.s. \quad $ so that $\quad \sup\limits_{u \in [0,t]} \left| S_u - \widetilde{S}_u^{I,\Gamma_n} \right| \underset{n\to \infty}{\to} 0 \ \PP - a.s.$. \myqed
\begin{remark}[Extension to semimartingales]\label{rk:semimartingale_extension}
\par \noindent In Theorems \ref{thm:p_partial_quantization_convergence} and \ref{thm:as_partial_quantization_convergence}, we limited ourselves to the case where $X$ is a local martingale. The proofs are easily extended to the case of a semimartingale $X$ as soon as there exists a locally finite measure $\nu$ on $[0,T]$ such that for every $\omega \in \Omega$ the finite-variation part $dV(\omega)$ in the canonical decomposition of $X$ is absolutely continuous with respect to $\nu$. In particular, this is the case for the Brownian bridge and Ornstein-Uhlenbeck processes whose finite-variation parts are absolutely continuous with respect to the Lebesgue measure on $[0,T]$.
\end{remark}
\begin{appendices}
\section{Injectivity properties of the Wiener integral}\label{sec:wiener_appendix}
\par In this appendix, we recall some results on the definition of the Wiener integral with respect to a Gaussian process. We focus on the injectivity properties. Here, we pay special attention to the special case of the Ornstein-Uhlenbeck processes.
\vspace{2mm}
\par \textbf{The covariance operator and the Cameron-Martin space}
\par \noindent Consider $X$ a bi-measurable centered Gaussian process on $[0,T]$ such that $\int_0^T \E[X_t^2]dt < \infty$ and with a continuous covariance function $\Gamma^X$ on $[0,T]\times [0,T]$. We denote by $H:=\overline{\sspan\left\{X_t, t\in [0,T] \right\}}^{L^2(\PP)}$ the Gaussian Hilbert space spanned by $(X_t)_{t \in [0,T]}$. The covariance operator $C_X$ of $X$ is defined by
$$
\begin{array}{lccc}
C_X :	& L^2([0,T])	& \to & L^2([0,T])\\
		& y				& \mapsto & C_X y = \E\left[(y,X) X\right].
\end{array}
$$
\par \noindent We have $C_X y(t) = \E\left[(y,X) X\right](t) = \E\left[\int_0^T X_s y(s) ds X_t\right] = \int_0^T \Gamma^X(t,s) y(s) ds$ where $\Gamma^X(t,s) = \E[X_tX_s]$ is the covariance function of $X$.
\par \noindent The Cameron-Martin space of $X$, (or reproducing Hilbert space of $C_X$), which we denote by $K_X$, is the subspace of $L^2([0,T])$ defined by $K_X := \left\{ t \mapsto \E\left[Z X_t \right], Z \in H\right\}$. $K_X$ is equipped with the scalar product defined by
$$
\langle k_1,k_2 \rangle_X = \E\left[Z_1 Z_2\right] \quad \textnormal{if} \quad k_i = \E\left[Z_i X_\cdot \right], \ i=1, 2,
$$
so that $(K_X,\langle \cdot \rangle_X)$ is a Hilbert space, isometric with the Hilbert space $\overline{\left\{(y,X) : y \in L^2([0,T])\right\}}^{H}$. $K_X$ is spanned as a Hilbert space by $\left\{C_X(y): y \in L^2([0,T]) \right\}$. 
\vspace{2mm}
\par \textbf{The Wiener integral}
\par \noindent Here, we follow the same steps as Lebovits and Lévy-Véhel in \cite{LebovitsLevyVehelCalSto} and Jost in \cite{JostVolterra} for the definition of a general Wiener integral. The difference here is that we use the quotient topology in order to define the Wiener integral in a more general setting. 
\par \noindent We define the map $U : H \to K_X$ defined by $U(Z)(t) = \E[Z X_t]$. By definition of $H$ and $K_X$, $U$ is a bijection and for any $s \in [0,T]$, we have $U(X_s) = \Gamma^X(s,\cdot)$. Consequently, $K_X$ is spanned by $\left(\Gamma^X(s,\cdot)\right)_{s \in [0,T]}$ as a Hilbert space. Now, we linearly map the set of the piecewise constant functions $\mathcal{E}([0,T])$ to the Cameron-Martin space $K_X$ by
$$
\begin{array}{lccc}
J:	& \mathcal{E}([0,T])	& \to & K_X \\
	& \1_{|s,t|}			& \mapsto & \Gamma^X(t,\cdot) - \Gamma^X(s,\cdot),
\end{array}
$$
where $|a,b|$ stands either for the interval $[a,b]$, $(a,b)$, $(a,b]$ or $[a,b)$. We equip $\mathcal{E}([0,T])$ with the bilinear form $\langle \cdot, \cdot \rangle_J$ which is defined by 
$$
\left\langle f, g \right\rangle_J := \left\langle Jf,Jg \right\rangle_X.
$$
\par \noindent It is a bilinear symmetric positive-\emph{semidefinite} form.
\begin{remark}
\par \noindent The so-called reproducing property shows that $\left\langle \1_{|0,t|}, \1_{|0,s|} \right\rangle_J = \Gamma^X(t,s) + \Gamma^X(0,0) - \Gamma^X(0,s) - \Gamma^X(0,t)$. When $X_0 = 0 \ a.s.$, this gives $\left\langle \1_{|0,t|}, \1_{|0,s|} \right\rangle_J = \Gamma^X(s,t)$. 
\end{remark}
\par Now, we define the equivalence relation $\underset{J}{\sim}$ on $\mathcal{E}([0,T])$ by $x \underset{J}{\sim} y$ if $\left\langle x-y, x-y \right\rangle_J = 0$. On the quotient space $E([0,T]) := \mathcal{E}([0,T]) / \underset{J}{\sim}$, the bilinear form $\langle \cdot, \cdot \rangle_J$ is positive-definite and thus it is a scalar product on $E([0,T])$. In this context, $J$ defines an (isometric) linear map from $E([0,T])$ to $K_X$. Then, considering the completion $F$ of $E([0,T])$ associated with this scalar product, $J$ is extended to $F$ and $U^{-1} \circ J: F \to H$ is an (isometric) injective map that we call Wiener integral associated to $X$.
$$
\int_0^T f(t) dX_t := U^{-1} \circ J(f). 
$$
\par \textbf{Injectivity properties of the Wiener integral}
\par \noindent As we have just seen, the Wiener integral is an (isometric) injective map from $F$ to $H$. Still, for example, when dealing with a standard Brownian bridge on $[0,T]$, $\left\|\1_{[0,T]}\right\|_J = 0$, so that there are functions of $\mathcal{E}([0,T])$ which have a nonzero $L^2$ norm and a zero $\|\cdot \|_J$ norm. Injectivity only holds in the quotient space $E([0,T]) = \mathcal{E}([0,T]) / \underset{J}{\sim}$ and its completion $F$. 
\vspace{2mm}
\par \noindent It is classical background that in the special case of a standard Brownian motion, $\|\cdot \|_J$ exactly coincides with the canonical $L^2$ norm so that $F = L^2([0,T])$. 
\vspace{2mm}
\par \textbf{Study of the case of Ornstein-Uhlenbeck processes}
\par \noindent From now, we will assume that $X$ is a centered Ornstein-Uhlenbeck process defined on $[0,T]$ by the SDE
$$
dX_t = -\theta X_t dt + \sigma dW_t \quad \textnormal{with } \sigma > 0 \textnormal{ and } \theta >0,
$$
where $W$ is a standard Brownian motion and $X_0 \simdist \mathcal{N}\left(0,\sigma_0^2\right)$ is independent of $W$. We make the additional assumption that $\theta T \leq \frac{4}{3}$. The covariance function writes
$$
\Gamma^X(s,t) = \frac{\sigma^2}{2\theta}e^{-\theta(s+t)} \left( e^{2 \min(s,t)}-1\right) + \sigma_0^2 e^{-\theta (s+t)}.
$$
\begin{prop}[Semi-norm equivalence on $\mathcal{E}({[}0,T{]})$]\label{prop:seminorm_domination}
\par There exist two positive constants $c$ and $C$ such that for every $f \in \mathcal{E}([0,T])$, $c\| f \|_2 \leq \| f \|_J \leq C \| f \|_2$.
\end{prop}
\par \noindent \textbf{Proof:} Let us consider $f \in \mathcal{E}([0,T])$. We have
$$
\|f\|_J^2	= \var\left( -\theta \int_0^T f(s) X_s ds + \sigma\int_0^T f(s) dW_s\right) \leq 2 \var\left( \theta \int_0^T f(s) X_s ds\right) + 2 \var\left( \sigma\int_0^T f(s) dW_s\right).
$$
The solution of the Ornstein-Uhlenbeck SDE is $X_t = X_0 e^{-\theta t} \overset{\independent }{+} \underbrace{\int_0^t \sigma e^{\theta (s-t)} dW_s}_{=: X_t^0}$. The so-defined process $\left(X_t^0\right)_{t \in [0,T]}$ is a centered Ornstein-Uhlenbeck process starting from $0$. Hence we have
$$
\begin{array}{lll}
\|f\|_J^2	&\leq 2 \var\left( X_0 \theta \int_0^T f(s) e^{-\theta s} ds\right) + 2\var\left( \theta \int_0^T f(s) X^0_s ds\right) + 2 \var\left( \sigma\int_0^T f(s) dW_s\right)\\
			&\leq 2 \theta^2 T \var(X_0) \int_0^T f(s)^2 ds + 2\var\left( \theta \int_0^T f(s) X^0_s ds\right) + 2 \var\left( \sigma\int_0^T f(s) dW_s\right).
\end{array}
$$
\par \noindent As in the proof of Proposition \ref{prop:OU_injectivity}, using that $\theta T < 4/3$, we can show that $\var\left( \theta \int_0^T f(s) X^0_s ds\right) \leq \var\left( \sigma\int_0^T f(s) dW_s\right)$. Hence
$$
\|f\|_J^2 \leq \underbrace{\left(2 \theta^2 T \sigma_0^2 + 4\sigma^2\right)}_{=: C^2} \int_0^T f(s)^2 ds,
$$
\par \noindent which is the desired inequality. Now we write
$$
\int_0^t f(s) dX_s = \underbrace{-\theta\int_0^T f(s) X_0 e^{-\theta s} ds}_{=: G^f_0} + \underbrace{\left(- \theta\int_0^T f(s) X_s^0 ds \right)}_{=: G^f_1} + \underbrace{\sigma \int_0^T f(s) dW_s}_{=: G^f_2},
$$
where $\left(G^f_0,G^f_1,G^f_2\right)$ is Gaussian and $G^f_0$ is independent of $G^f_1$ and $G^f_2$. Hence 
\begin{multline}\label{eq:intermediate_seminorm_domination}
\var\left(\int_0^t f(s) dX_s \right) \geq \var\left(G^f_1 + G^f_2\right) = \var\left(G^f_1\right)+\var\left(G^f_2\right) + 2\cov\left(G^f_1,G^f_2\right)\\
								\geq \var\left(G^f_1\right)+\var\left(G^f_2\right) -2\sqrt{\var\left(G^f_1\right)\var\left(G^f_2\right)} = \left(\sqrt{\var\left(G^f_2\right)}-\sqrt{\var\left(G^f_1\right)}\right)^2.
\end{multline}
\par It has been shown at the beginning of the proof of Proposition \ref{prop:OU_injectivity} that there exists a constant $K<1$ independent of $f$ such that
$\var\left(G^f_1\right) \leq K \var\left(G^f_2\right)$. $K$ was defined by
$$
K = \frac{\theta^2}{\sigma^2} \sqrt{\int_0^T \int_0^T \left(\Gamma^{X^0}(s,t)\right)^2 ds dt},
$$
where $\Gamma^{X^0}$ is the covariance function of the Ornstein-Uhlenbeck process starting from $0$. Plugging this into Equation (\ref{eq:intermediate_seminorm_domination}) yields
$$
\var\left(\int_0^t f(s) dX_s \right) \geq \left(1-\sqrt{K}\right)^2 \var\left(G^f_2\right) = \underbrace{\left(1-\sqrt{K}\right)^2 \sigma^2}_{=: c^2} \|f\|_2^2.
$$
This is the wanted inequality. \myqed
\vspace{2mm}
\par A straightforward consequence of Proposition \ref{prop:seminorm_domination} is that $\|f\|_J = 0 \Leftrightarrow \|f\|_2=0$ so that equivalence classes in $\mathcal{E}([0,T])$ for the relation $\underset{J}{\sim}$ are \emph{almost surely} equal functions. Another consequence is that the sets of Cauchy sequences and convergent sequences for the two norms on $E([0,T])$ coincide, and thus the corresponding completions of $E([0,T])$ are the same. In other words, in the case of Ornstein-Uhlenbeck processes that satisfy the condition $\theta T \leq \frac{4}{3}$, we have $F = L^2([0,T])$. 
\end{appendices}
\vspace{6mm}
\par \noindent The author is grateful to Benjamin Jourdain and Gilles Pagès for their helpful remarks and comments.
\bibliography{biblio}
\end{document}